\documentclass[11pt,a4paper]{article}
\usepackage[T1]{fontenc}
\usepackage[utf8]{inputenc}
\usepackage{lmodern}
\usepackage[margin=2.45cm]{geometry}
\usepackage{amsmath,amssymb,amsthm}
\usepackage{graphicx,booktabs,enumitem}
\usepackage{microtype}
\usepackage{xcolor}
\usepackage[colorlinks=true,linkcolor=blue!40!black,citecolor=blue!40!black,urlcolor=blue!40!black]{hyperref}
\hypersetup{pdftitle={On the V-toothpick sequence},pdfauthor={Benoit Cloitre}}
\newtheorem{theorem}{Theorem}
\newtheorem{lemma}[theorem]{Lemma}
\newtheorem{proposition}[theorem]{Proposition}
\newtheorem{corollary}[theorem]{Corollary}
\theoremstyle{definition}
\newtheorem{definition}[theorem]{Definition}

\theoremstyle{remark}
\newtheorem{remark}[theorem]{Remark}
\newcommand{\lb}{\operatorname{lb}}
\newcommand{\NN}{\mathbb N}
\newcommand{\ZZ}{\mathbb Z}

\newcommand{\bitand}{\mathbin{\&}}
\setlist{itemsep=2pt,topsep=4pt}
\title{On the V-toothpick sequence}
\author{Beno\^it Cloitre\\[3pt]\small\href{https://orcid.org/0009-0001-6778-153X}{ORCID 0009-0001-6778-153X}}
\date{}
\begin{document}
\maketitle
\vspace{-1.3em}
\begin{abstract}
We give a generating function for the V-toothpick sequence and an
algorithm that computes any term in time polynomial in the number of binary
digits of its index. We determine its quadratic growth and periodic
fluctuations, and give convergent series for the constants in its asymptotic bounds.
\end{abstract}

\section*{Notation and conventions}

We write $\NN=\{0,1,2,\ldots\}$ and $\ZZ$ for the integers, and use $\log_2$ for the
binary logarithm and $\log$ for the natural logarithm. An integer is \emph{dyadic} if it
is a power of two. For $m\ge1$, $\lb(m)$ is the largest power of two dividing $m$.
For $m\ge0$, $s_2(m)$ counts the ones in its binary expansion. For nonnegative integers,
$x\bitand y$ has binary digits equal to the products of the corresponding digits of $x$ and $y$.
Thus $x\bitand y=0$ when they have no common digit equal to one.

Generating-function identities are understood as identities of formal power series.
We write $[z^j]Q$ for the coefficient of $z^j$ in $Q$, with value zero for $j<0$.
For a polynomial $Q$, $\operatorname{supp}(Q)$ is the set of its nonzero coefficient indices,
$\#\operatorname{supp}(Q)$ its size, and $\|Q\|$ the sum of the absolute values of its coefficients.
The indicator $\mathbf1_{\mathcal E}$ is $1$ when the condition $\mathcal E$ holds and $0$ otherwise.
For $n\ge0$, we set $\binom nj=0$ unless $0\le j\le n$.

\section{Introduction}

\subsection{From toothpicks to the honeycomb}

Omar E.~Pol introduced the toothpick sequence, OEIS A139250, in 2008 \cite{OEIS}.
Starting with one segment on the square lattice, the construction adds a perpendicular segment whose
midpoint meets each exposed end. Applegate, Pol, and Sloane \cite{APS} studied this sequence in 2010,
giving recurrences based on powers of two, generating functions, and asymptotic results. They also related the growth to
the Ulam-Warburton cellular automaton and examined other local rules.

Pol introduced the V-toothpick variant in 2009 as A161206 \cite{OEIS}. Its pieces lie on the honeycomb
lattice, and the counting mechanism of the square-lattice construction does not directly apply.
The OEIS entry describes the growth but gives no general formula for its terms. We obtain an expression
for its generating function, an algorithm for an arbitrary term, and its leading asymptotic form.
The algorithm computes a term without computing all preceding terms. To our knowledge, this is the first
such representation and algorithm with polynomial cost in the binary input length for this variant.

Consider the honeycomb lattice, that is, the tiling of the plane by regular hexagons with sides of length
one. A \emph{V-toothpick} is a piece made of two sides of the lattice that meet at a vertex at an angle of
120 degrees. This vertex is the apex of the piece and the two other endpoints are its ends.

\begin{figure}[ht]
\centering
\includegraphics[width=\linewidth]{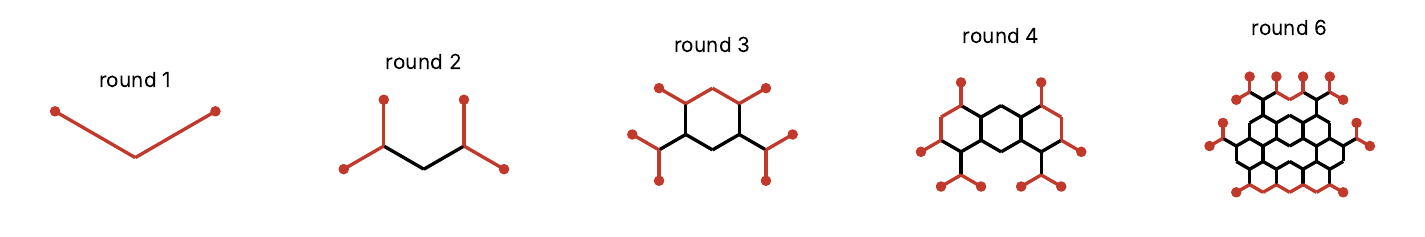}
\caption{The structure after rounds 1, 2, 3, 4, and 6. The pieces of the last round are in red and the dots
are the free ends.}
\label{fig:rounds}
\end{figure}

\begin{figure}[ht]
\centering
\includegraphics[width=0.62\linewidth]{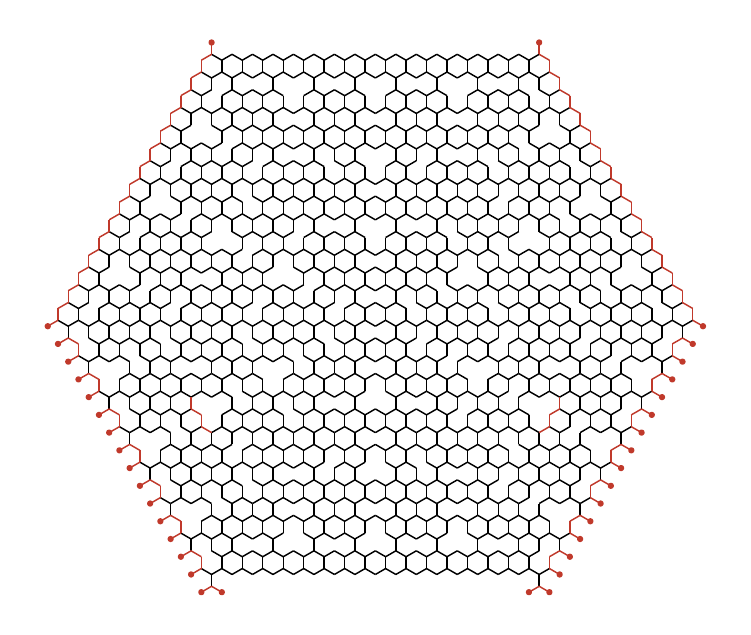}
\caption{The V-toothpick structure after 32 rounds. It contains $a(32)=1005$ pieces.}
\label{fig:structure}
\end{figure}

The game starts with an empty plane. At round~1 one piece is placed. At each later round a new piece is
placed at every free end of the structure. A \emph{free end} is an end of a piece that touches no other
piece. The new piece has its apex at the free end and covers the two other sides of the lattice that leave
this vertex. All pieces of a round are placed simultaneously. When two new pieces reach the same vertex
in the same round, this vertex is not a free end and nothing grows from it later. Figure~\ref{fig:rounds}
shows the first rounds and Figure~\ref{fig:structure} shows the structure after 32 rounds.

Let $a(n)$ be the number of pieces after $n$ rounds, with $a(0)=0$, and let $t(n)=a(n)-a(n-1)$ be the number
of pieces placed at round $n$. The first values are
\[
a(n)=0,\,1,\,3,\,7,\,13,\,21,\,31,\,43,\,57,\,69,\,81,\,99,\,123,\ldots
\]
The increments form entry A161207 of the OEIS \cite{OEIS}.

\subsection{The question}

The definition gives $a(n)$ by playing the game through all the preceding
rounds. Finite propagation places every affected vertex in a region with
$O(n^2)$ sites. A direct simulation visits successive fronts in this region,
whereas the question here is whether one can recover the count from the
binary expansion of $n$. Even the number of simulated rounds is exponential
in the length of that input.

\begin{figure}[ht]
\centering
\includegraphics[width=0.8\linewidth]{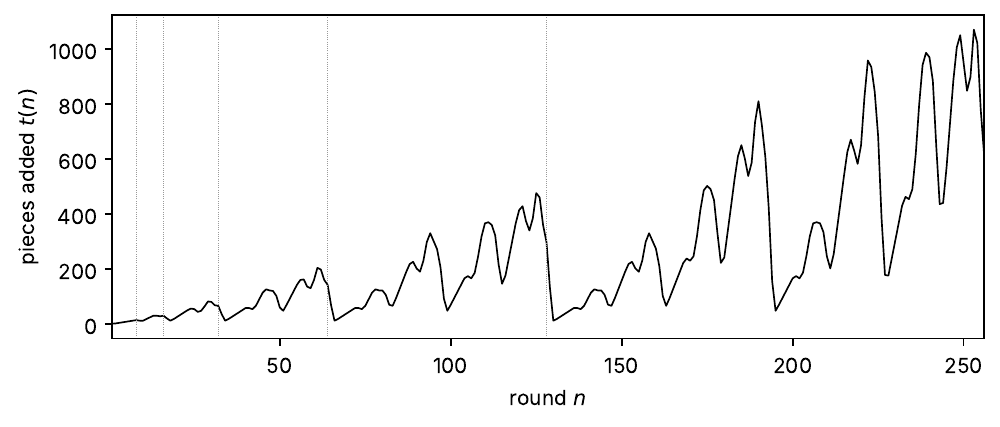}
\caption{The number $t(n)$ of pieces placed at round $n$, for $n\le 256$. The dotted lines mark the powers of
two.}
\label{fig:increments}
\end{figure}

The irregular increments in Figure~\ref{fig:increments} suggest that
the dyadic repetition is more complicated than for the original toothpick
sequence. We seek a description that accounts for all vertices that appear inside
each region, including those reached when growing parts of the structure meet.

\subsection{Results}

We use polynomials to count the vertices that appear at successive times under the local growth rule
of Section~\ref{sec:statements}. Each appearance is called a \emph{birth}.
The polynomials $P_2,P_4,P_8,\ldots$ are indexed by the powers of two. The polynomial $P_q$ counts, round by round, the births in
a finite region of the structure that we call a diamond of size $q$. Its definition is given in
Section~\ref{sec:statements}. The first two are $P_2=1+z$ and $P_4=1+2z+2z^2+z^3+2z^4+2z^5$.
The cone polynomial $C_r$ counts births in a triangular region of size $r$,
also defined there. Set
\[
 \mathcal G(w)=\prod_{i\ge0}(1+2w^{2^i})
              =\sum_{j\ge0}2^{s_2(j)}w^j.
\]

\begin{theorem}\label{thm:series}
Let $A(z)=\sum_{n\ge 0}a(n)z^n$. In the ring of formal power series,
\begin{equation}\label{eq:SG}
A(z)=\frac{(1+z)(1+2z^2)}{1-z}\left(\frac{z}{1-z}+2\sum_{q=2,4,8,\ldots}\frac{z^{2q+1}P_q(z)}{1-z^{4q}}\right).
\end{equation}
The equivalent expression in terms of cones is
\[
 A(z)=\frac{(1+z)(1+2z^2)}{1-z}
 \left(\frac{z}{1-z}+2\sum_{q=2,4,8,\ldots}
 z^{2q}(C_q(z)-1)\,\mathcal G(z^{4q})\right).
\]
\end{theorem}

Each term of the first sum says that a diamond of size $q$ appears for the first time at a fixed round and then
again every $4q$ rounds. Only the terms with $2q+1\le n$ contribute to the coefficient of $z^n$, so the sum
is well defined.

For example, the contribution of $q=2$ to the increment series $(1-z)A(z)$ begins with
\[
 2z^5(1+z)^2(1+2z^2)=2z^5+4z^6+6z^7+8z^8+4z^9.
\]
The rational term contributes $6$ at each round $n\ge4$. Hence $t(5)=8$, $t(6)=10$, and $t(7)=12$,
which give $a(5)=21$, $a(6)=31$, and $a(7)=43$.

The geometric construction in Section~\ref{sec:geom} shows that $P_q$ has at least
$2q-2$ nonzero coefficients. Expanding a polynomial of size comparable to $n$ would therefore
require at least as many operations. The following algorithm avoids this expansion.
The class FP consists of functions computable in time polynomial in the length of the input~\cite{AB}.

\begin{theorem}\label{thm:algo}
There is an algorithm that takes as input the binary expansion of an integer $n\ge 0$ and returns $a(n)$ and
$t(n)$. With $\ell=\lfloor\log_2(n+3)\rfloor+1$, it performs $O(\ell^8)$ arithmetic operations on integers of
$O(\ell)$ bits. In particular the maps $n\mapsto a(n)$ and $n\mapsto t(n)$ belong to the class FP.
\end{theorem}

Here the input has length of order $\log n$, so the running time is polynomial in $\log n$. The
same program works for every $n$ and needs no precomputed table. We do not claim that the exponent 8 is optimal. The computations in Section~\ref{sec:exp} include indices up to $10^{100}$,
with the outputs checked by two separately written implementations.

The number of births in a diamond of size $q$ is
asymptotic to $\kappa q^2$ for a computable constant $\kappa=0.6879156855\ldots$. Theorem~\ref{thm:growth}
shows that
\[
 a(n)=n^2\,\Psi\Bigl(\frac{n}{2^k}\Bigr)+O\bigl(n^{\log_23}(\log n)^2\bigr),\qquad k=\lfloor\log_2n\rfloor,
\]
for a continuous function $\Psi$ on $[1,2]$ determined by two functional equations. The ratio $a(n)/n^2$ has
no limit. Its lower and upper limits are $9\kappa/8$ and $3\kappa/2$.
Theorem~\ref{thm:ternary} gives a series for $\kappa$ whose first $N$ terms
have error at most $4/3^N$ and require cones of size at most $4N$.

Theorem~\ref{thm:windows} also shows that the differences $X_k=C_{2^k}-C_{2^k-1}$ have only
$O(k^3)$ nonzero coefficients. It computes these coefficients using cones and triangles of size $O(k)$.

The recursive construction computes its initial data by applying the growth rule to triangles of side
$O(\log n)$. Whether a fixed finite family of numerical recurrences suffices remains open.

\subsection{Geometric reduction and coefficient extraction}

Earliest births follow Pascal's triangle modulo $2$. Along certain lines of constant coordinate sum,
this pattern prevents further growth and separates the structure into finite diamonds. The repeated occurrence of a
diamond of size $q$, every $4q$ rounds, gives the denominator in \eqref{eq:SG}.

Each diamond reduces to smaller diamonds and cones. A translation identity for truncated triangles
then reduces the cones, except for domains of size $O(\log n)$. The algorithm stores these identities
as a graph of polynomial expressions. A coefficient or a partial sum is evaluated by following this
graph. The queried indices lie in a polynomial number of residue classes modulo the dyadic sizes,
which bounds the number of evaluations.

\subsection{Related work}

The square-lattice analysis of Applegate, Pol, and Sloane \cite{APS} provides
the closest enumerative comparison. Here the relevant finite regions also
occur at dyadic scales, but their evolution requires additional cone
polynomials. Their reduction gives a uniform algorithm as the queried index grows.

The complexity of prediction for cellular automata depends on the rule and the input model.
Neary and Woods~\cite{NW} show that predicting Rule~110 can encode any decision problem solvable
in polynomial time, through a transformation using logarithmic space. Such results concern families
of initial configurations and a prediction task. Our input consists only
of a time index, with a fixed initial seed, and our output is a global
count. We make no claim about prediction from arbitrary initial data.

The representation used below is related to random access in compressed
objects. Bille et al.\ \cite{Bille} retrieve characters from a string specified by a grammar
that generates that string alone, without expanding it. Our nodes instead describe signed
sums of shifted, potentially overlapping polynomials. The residue lemma
is needed to bound how many coefficient queries these sums generate.

\section{The model and the birth polynomials}\label{sec:statements}

\subsection{The lattice and the growth rule}

All time indices start at $0$, unless global round numbers are specified.
We describe the honeycomb lattice as a bipartite graph. Its vertices are the points $A(x,y)$ and $B(x,y)$
with $(x,y)\in\ZZ^2$. We also write $A(x,y)=(0,x,y)$ and $B(x,y)=(1,x,y)$, so that a vertex is a triple
$(p,x,y)$ with $p\in\{0,1\}$. The neighbors of $A(x,y)$ are $B(x,y)$, $B(x-1,y)$, and $B(x,y-1)$.
Consequently the neighbors of $B(x,y)$ are $A(x,y)$, $A(x+1,y)$, and $A(x,y+1)$.

\begin{definition}\label{def:growth}
Let $D$ be a set of vertices and let $s\in D$ be a vertex called the seed. Time is discrete. Every vertex of
$D$ other than the seed is initially \emph{intact}. The seed is born at time $0$. At each time $t\ge 0$ every
vertex born at time $t$ sends a proposal to each of its neighbors in $D$. An intact vertex that receives
exactly one proposal at time $t+1$ is \emph{born} at time $t+1$. An intact vertex that receives two or three
proposals at time $t+1$ is \emph{blocked}. A vertex that is born or blocked is never modified again and
ignores later proposals. The \emph{birth polynomial} of $(D,s)$ is the sum of $z^{t}$ over the vertices born
at time $t$. When $D$ is infinite it is a formal power series.
\end{definition}

A born vertex emits once, at the time that follows its birth. A blocked vertex never emits.

\begin{remark}\label{rem:finite-domain}
In a finite domain with $N$ vertices, a birth at time $t\ge1$ requires a birth at time $t-1$.
Thus the birth polynomial has degree at most $N-1$.
\end{remark}

\begin{lemma}\label{lem:pieces}
Run the growth of Definition~\ref{def:growth} on the whole lattice with seed $A(0,0)$, with the single change
that the seed sends proposals to two of its three neighbors only. Then for every $n\ge 1$ the number $t(n)$
of pieces placed at round $n$ is the number of vertices born at time $n-1$.
\end{lemma}

\begin{proof}
Place the first piece with its apex at $A(0,0)$ and its two ends at the two chosen neighbors. We prove by
induction on $n\ge 1$ that after round $n$ the vertices touched by the structure are the vertices that are
not intact at time $n$, and that the free ends are the vertices born at time $n$. After round $1$ the two
ends of the first piece have each received one proposal and are born at time $1$, and they are the free
ends. Assume the claim after round $n$. At round $n+1$ a piece is placed at each free end $u$. It covers the
two sides at $u$ that are not yet covered, which lead to the two neighbors of $u$ other than the vertex that
sent its proposal to $u$. In the growth, $u$ sends a proposal to its three neighbors and the third one is
ignored because its target is not intact. An end $w$ of a new piece is a free end after round $n+1$ exactly
when $w$ was not touched before and is reached by a single new piece. This means that $w$ is intact at time
$n$ and receives exactly one proposal at time $n+1$, that is, $w$ is born at time $n+1$. The vertices touched
for the first time at round $n+1$ are those that receive their first proposals at time $n+1$. This proves the
claim. The pieces placed at round $n+1$ are in bijection with the free ends after round $n$, and the piece of
round~$1$ corresponds to the seed.
\end{proof}

\subsection{Cones and derived polynomials}

For $r\ge 1$ the \emph{cone} of size $r$ is
\[
K_r=\{(p,x,y)\ \mid\ x\ge 0,\ y\ge 0,\ 2x+y+p<r\}.
\]
It has $r(r+1)/2$ vertices. We denote by $C_r(z)$ the birth polynomial of the growth on $K_r$ with seed
$A(0,0)$, and we set $C_0=0$. For instance
\[
C_1=1,\qquad C_2=1+z+z^2,\qquad C_3=1+z+2z^2+z^3+z^4,\qquad C_4=1+z+2z^2+2z^3+z^4+z^5+z^6 .
\]
For integers $u\le v$ we write $J_{u,v}=z^u+z^{u+1}+\cdots+z^v$, and $J_{u,v}=0$ when $u>v$. We set
\[
H_r=C_r+zC_{r-1}-J_{1,r-1},
\]
with $H_0=0$.

\begin{definition}\label{def:P}
For every dyadic $q\ge2$, define $P_q$ recursively by
\begin{equation}\label{eq:diamond-cone}
 P_q=S_{q,0}+\frac{C_q-1}{z},
\end{equation}
where, for dyadic $B\ge1$ and an integer $d\ge0$,
\begin{equation}\label{eq:S}
 S_{B,d}=\sum_{\substack{2\le m<B\\m\ \mathrm{even}}}z^{2m+d}P_{\lb(m)}.
\end{equation}
\end{definition}

The quotient is a polynomial because $C_q(0)=1$. The definition is well founded, since the sum contains
only $P_p$ with $p<q$. In particular $P_2=(C_2-1)/z=1+z$.
Proposition~\ref{prop:diamonds} identifies $P_q$ with the births in a finite diamond.
Lemma~\ref{lem:degrees} bounds its degree by $2q+O((\log q)^2)$.
The equality $\deg P_q=2q-3$ is proved over a finite range in Appendix~\ref{app:bounds} and remains open in general.

\subsection{Generating functions and computation model}

Theorem~\ref{thm:series} can now be read with Definition~\ref{def:P}. We also record the form for the
increments. Let $B(z)=\sum_{n\ge 1}t(n)z^n=(1-z)A(z)$ and
\begin{equation}\label{eq:G}
G(z)=\sum_{q=2,4,8,\ldots}\frac{z^{2q+1}P_q(z)}{1-z^{4q}} .
\end{equation}
Then \eqref{eq:SG} is equivalent to
\begin{equation}\label{eq:SGprime}
B(z)=(1+z)(1+2z^2)\left(\frac{z}{1-z}+2G(z)\right).
\end{equation}
Theorem~\ref{thm:algo} refers to the following model of computation. The input is the binary expansion of
$n$. An arithmetic operation is an addition, a subtraction, a multiplication, a comparison, or a division with
remainder on integers of $O(\ell)$ bits. The algorithm is uniform, which means that a single program handles
every $n$, builds its own data, and reads no table.

\section{From geometric barriers to the generating function}\label{sec:geom}

\subsection{Earliest times and Pascal's triangle}

For a vertex $(p,x,y)$ with $x,y\ge 0$ let $t_{\min}(p,x,y)=2(x+y)+p$. Along every edge of the lattice the value of
$t_{\min}$ changes by exactly one, and $t_{\min}$ vanishes at $A(0,0)$. In a growth seeded at $A(0,0)$ inside the quadrant
$x,y\ge 0$, no vertex can therefore receive a proposal before time $t_{\min}$. We call $t_{\min}$ the \emph{earliest time}
of the vertex.

We say that a domain $D$ inside the quadrant is \emph{closed under predecessors} when, for every vertex of $D$,
the neighbors with nonnegative coordinates and with a smaller value of $t_{\min}$ also belong to $D$. The cones
$K_r$ and the truncated triangles of Section~\ref{sec:cones} have this property, because a step that
decreases $t_{\min}$ increases neither $x+y$ nor $2x+y+p$.

\begin{lemma}\label{lem:pascal}
Let $D$ be a domain inside the quadrant that contains $A(0,0)$ and is closed under predecessors. In the growth
on $D$ seeded at $A(0,0)$, a vertex $A(x,y)$ or $B(x,y)$ of $D$ is born at its earliest time if and only if
$x\bitand y=0$.
\end{lemma}

\begin{proof}
A vertex can be born at its earliest time only through a proposal sent by a neighbor that was itself born at
its own earliest time, one time step earlier. For $B(x,y)$ the only neighbor with a smaller value of $t_{\min}$ is
$A(x,y)$, so $B(x,y)$ is born at its earliest time exactly when $A(x,y)$ is. For $A(x,y)$ the neighbors with
a smaller value of $t_{\min}$ are $B(x-1,y)$ and $B(x,y-1)$, when their coordinates are nonnegative, and they
belong to $D$ by assumption. Let $\alpha(x,y)$ be $1$ when $A(x,y)$ is born at time $2(x+y)$ and $0$
otherwise, with $\alpha=0$ outside the quadrant. The rule of a single first arrival gives
$\alpha(x,y)=\alpha(x-1,y)+\alpha(x,y-1)$ modulo~$2$ for $(x,y)\ne(0,0)$, since two simultaneous arrivals
block the vertex and no arrival leaves it intact at that time. With $\alpha(0,0)=1$ this is the recurrence of
the binomial coefficients, so $\alpha(x,y)$ is the parity of $\binom{x+y}{x}$. By the theorem of Lucas this
binomial coefficient is odd exactly when the addition of $x$ and $y$ in base two has no carry, which is the
condition $x\bitand y=0$.
\end{proof}

\subsection{Spines and a global counting identity}

We use round numbers for the global growth, so that the seed $A(0,0)$ is born at
round $1$ and its two chosen neighbors $B(0,0)$ and $B(-1,0)$ at round $2$.
Birth times in a finite cone are still measured from its seed at time $0$.
Put
\[
 g=2x+y+p,\qquad h=y-x,\qquad
 Q_1=\{x\ge0,\ y\ge0\},\qquad
 H^+=\{x\ge1,\ y>x\}.
\]
The vertical spine consists of $A(0,m)$ and $B(0,m)$ for $m\ge1$.
Let $Z(n)$ count the births at round $n$ in $\{x\ge0,y>x\}$ and let $F(n)$
count those in $H^+$. We put $F(n)=0$ for negative indices. We have
$Z(n)=1+F(n)$ for $n\ge3$.

\begin{lemma}\label{lem:spines}
The quadrant $Q_1$ evolves as the growth on the induced quadrant with seed
$A(0,0)$ at round $1$. Its vertical and horizontal spines have birth times
\[
 A(0,m),A(m,0)\colon 2m+1,\qquad
 B(0,m),B(m,0)\colon 2m+2.
\]
The ray $x=y>0$ is inactive at every round, as is the rest of its diagonal
inside $Q_1$ after round $2$. The global axis $g=0$ is inactive after the seed.
There is also a spine
\[
 B(m,-m)\colon 2m+2,\qquad A(m+1,-m)\colon 2m+3
 \quad(m\ge0).
\]
\end{lemma}

\begin{proof}
The functions $2x+p$, $2y+p$, and $2(x+y)+p$ change by one along each edge.
Their values at the two births of round $2$ give the lower bounds
\[
 n\ge2x+p+1,\quad n\ge2y+p+1,\quad n\ge2(x+y)+p+1
\]
for any later first arrival. Except at the already resolved initial seed,
the same argument with $-2y-p$ gives $n\ge-2y-p+3$.
On the vertical spine, each indicated birth has exactly one predecessor
at the required round. The other two neighbors have a lower bound two
rounds later than that predecessor. Induction proves the displayed times,
starting with $B(0,0)$, and the argument with $x$ and $y$ interchanged proves
the horizontal times. The only external edges of $Q_1$ away from the seed
join $B(-1,m)$ to $A(0,m)$ or $B(m,-1)$ to $A(m,0)$.
Their external endpoints cannot be born before round $2m+2$, so every return
arrives after the spine has been resolved. This proves autonomy.

The quadrant and its initial state are invariant under
$\mu(p,x,y)=(p,y,x)$. After $B(0,0)$, a diagonal vertex receives its
off-diagonal proposals in pairs. Its possible diagonal neighbor is
inactive by induction. Thus no further diagonal birth occurs. Likewise,
the global reflection
\[
 R(p,x,y)=(p,-x-y-p,y)
\]
preserves the initial V, exchanges its two ends, and sends $g$ to $-g$.
On its fixed axis the off-axis proposals arrive in pairs. The only edge
from the seed along that axis is the edge not used by the initial V.
Induction therefore proves inactivity of the axis after the seed.

The vertex $B(m,-m)$ reaches $A(m+1,-m)$ at round $2m+3$.
The other possible parents have lower bound $2m+4$. That A vertex then
reaches $B(m+1,-m-1)$ at round $2m+4$, whereas its other possible parents
have lower bound $2m+5$. These bounds also exclude any earlier first touch.
Induction from $B(0,0)$ proves the last spine.
\end{proof}

In particular, $Q_1$ has $2Z(n)$ births at round $n\ge3$, and
$Z(n)=1+F(n)$. Its earliest births are given by Lemma~\ref{lem:pascal},
with the additional shift of one round. Along the boundary of $H^+$,
$A(1,m)$ is born at round $2m+3$ for even $m\ge2$. For odd $m$ it is
blocked at that round by $B(0,m)$ and $B(1,m-1)$.
We keep these literal boundary proposals when using symmetry arguments.

\begin{proposition}
For every $n\ge4$,
\begin{equation}\label{eq:filter}
 t(n)=6+2F(n+2)+2F(n+1)+4F(n)+4F(n-1).
\end{equation}
\end{proposition}

\begin{proof}
By the reflection $R$, it suffices to count the births in $g>0$.
Besides $Q_1$, introduce
\[
 S^-=\{g>0,h>0,x\le-1\},\quad
 Q_2=\{y\le0,x+y+p\ge1\},\quad
 D=\{g>0,x+y+p\le0\}.
\]
These regions give a partition into $Q_1$, $S^-$, $Q_2\setminus Q_1$, and $D$
of the births in $g>0$ after round $3$.

The graph automorphism
\[
 \Phi(p,x,y)=(1-p,-x,x+y+p)
\]
maps $S^-$ onto $H^+$ and sends $(g,h)$ to $(h+1,g)$.
The proposals entering $S^-$ are those from $A(0,m)$ to $B(-1,m)$.
Their images are the proposals from $B(0,m)$ to $A(1,m)$ entering $H^+$,
one round later. In both cases $m\ge2$ for an interior target.
The other boundary is an inactive ray, and returns to the spine arrive
after its resolution by Lemma~\ref{lem:spines}. The two regions are empty
off their boundaries at corresponding rounds $3$ and $4$. Induction on
first arrivals, including simultaneous arrivals, proves that $S^-$ has
$F(n+1)=Z(n+1)-1$ births at round $n$.

The region $Q_2$ is autonomous. Its external edges join $B(x,0)$ to
$A(x,1)$, or $B(m,-m)$ to $A(m,-m)$. On the first boundary the spine
vertex is born at $2x+2$ and the external A cannot be born before $2x+3$.
On the second boundary the corresponding bounds are $2m+2$ and $2m+3$,
with the initial seed already resolved when $m=0$.
Every return is therefore too late. The automorphism
\[
 K(p,x,y)=(1-p,x+y+p,-x)
\]
maps $Q_1$ onto $Q_2$. It maps the state of $Q_1$ at round $2$ to that
of $Q_2$ at round $3$, so that $Q_2$ has $2Z(n-1)$ births.
Its intersection with $Q_1$ is the horizontal spine and has one birth
per round. Thus $Q_2\setminus Q_1$ contributes $2Z(n-1)-1$.

The remaining automorphism
\[
 J(p,x,y)=(p,1-x-y-p,x+1)
\]
maps $D$ onto $H^+$. Its new coordinates satisfy $x'=1-(x+y+p)$
and $h'=g$. It takes the inactive boundary $g=0$ to $h'=0$ and maps
the entering edge $B(m,-m)A(m,-m)$ to
$B(0,m+1)A(1,m+1)$. The source of the first edge is born at $2m+2$
and the source of the image edge at $2m+4$.
These are all the entering edges for $m\ge1$. The otherwise missing
edge at $m=0$ has its image target on the inactive diagonal.
The regions are empty at corresponding rounds $3$ and $5$.
Their states therefore agree after this translation of time, and $D$
has $F(n+2)=Z(n+2)-1$ births at round $n$.

Adding the four counts gives
\[
 \frac{t(n)}2=Z(n+2)+Z(n+1)+2Z(n)+2Z(n-1)-3.
\]
Substitution of $Z=1+F$ proves the formula.
\end{proof}

\subsection{Isolation and coverage of the diamonds}

For an even integer $m\ge2$, put $q=\lb(m)$ and use the coordinates
$\ell=y-x$ and $r=x+y$. Define the finite region $D_m$ by
\begin{equation}\label{eq:diamond-domain}
 m-q<\ell<m\le r<m+q,
 \qquad\text{with A vertices on }r=m\text{ omitted}.
\end{equation}
Its seed is $A(1,m)$. After writing $w=y-m$, the region depends only on $q$
and becomes
\begin{equation}\label{eq:standard-diamond}
 x\ge1,\qquad -x\le w\le x-1,\qquad x+|w|\le q-1,
 \qquad A(x,-x)\text{ omitted}.
\end{equation}
Write $\widehat P_q$ for its birth polynomial when the seed is born at age $0$.

\begin{lemma}\label{lem:isolation}
The growth in $D_m$ is its isolated growth, with seed born at round $2m+3$.
Every birth in $H^+$ belongs to exactly one such region.
\end{lemma}

\begin{proof}
Two families of blocked vertices separate the diamonds. If $M$ is even and
$Q=\lb(M)$, then for $0<x<Q$ the vertex $A(x,M-x)$ has two earliest
parents $B(x-1,M-x)$ and $B(x,M-x-1)$ at round $2M$.
Their coordinates have disjoint binary digits. Lemma~\ref{lem:pascal}
and the lower bound on first arrivals therefore block the A vertex at
round $2M+1$, permanently.

For the second barrier fix even $m$, put $q=\lb(m)$, and translate
$y=m+w$. The triangle $x,w\ge0$, $x+w<q$ is empty immediately before
round $2m+1$, when only $A(0,m)$ is born. This follows from the lower
bound $2m+2(x+w)+p+1$.
On its horizontal axis the vertices are born at their earliest times,
since $x\bitand m=0$ for $x<q$. A possible return from below comes from
$B(x,m-1)$. For $0<x<q$ one has $x\bitand(m-1)=x\ne0$.
That B misses its earliest time and, by bipartite parity, any later birth
sends its proposal after $A(x,m)$ has been resolved.
The vertical axis is the one-way spine already proved.

The far side $x+w=q$ is closed by the first barrier because
$\lb(m+q)\ge2q$. At its endpoints the only interior neighbors are
$B(0,m+q-1)$ and $B(q-1,m)$, both born at round $2(m+q)$.
The endpoint A vertices cannot be born before round $2(m+q)+1$, so their
returns are also ineffective. The induced triangle thus has a single
seed, symmetric initial state, and no effective external input.
Reflection in $x=w$ preserves its evolution. A positive diagonal
vertex receives its off-diagonal proposals in pairs, and its diagonal
neighbor is inactive by induction. This proves that
\[
 y-x=m,\qquad m<x+y<m+q
\]
is an inactive segment at every round. This argument uses the actual
initially intact sites, including those later blocked by spine proposals.

These barriers determine the four faces of \eqref{eq:diamond-domain}.
The faces $r=m$ and $r=m+q$ are closed by the A barriers just proved.
The face $\ell=m$ is closed by the mirror segment. If $m=q$, the face
$\ell=m-q$ is the original inactive diagonal. Otherwise put $m_0=m-q$.
Then $\lb(m_0)\ge2q$, and the relevant fourth face lies on the mirror
segment already proved for $m_0$. This uses no previous diamond decomposition.
At the corners, the remaining incident A vertices lie on these barriers.
The B vertices on $r=m$ have their same-level A neighbor blocked and
their other neighbors on $r=m+1$. On the far face no extra B vertex
has an interior neighbor. The only effective entering edge is therefore
$B(0,m)A(1,m)$, which creates the seed at round $2m+3$.
No other birth or first touch can occur inside before this input, so
the entire subsequent evolution is the isolated one.

For coverage, let a vertex of $H^+$ have coordinates $0<\ell<r$.
The integers $\ell$ and $r$ have the same parity. Let
\[
 q=2^{\lfloor\log_2(\ell\mathbin{\mathrm{xor}}r)\rfloor},
 \qquad m=q\lfloor r/q\rfloor,
\]
where xor denotes digitwise exclusive or. The highest differing binary
digit is $0$ in $\ell$ and $1$ in $r$. Thus $q\ge2$, $\lb(m)=q$ and
$m-q\le\ell<m\le r<m+q$. Equality at $\ell=m-q$, or an A vertex at
$r=m$, places the site on one of the inactive barriers. Otherwise it
belongs to $D_m$. The highest differing digit also makes this choice
unique. Hence the isolated diamonds contain all births and are disjoint.
\end{proof}

\subsection{The polynomial of a diamond}

\begin{proposition}\label{prop:diamonds}
The geometric polynomial $\widehat P_q$ is $P_q$ of Definition~\ref{def:P}.
Consequently
\begin{equation}\label{eq:F}
 F(z)=\sum_{q=2,4,8,\ldots}\frac{z^{2q+3}P_q(z)}{1-z^{4q}}.
\end{equation}
\end{proposition}

\begin{proof}
Let $q\ge2$ be dyadic. In the domain
\eqref{eq:standard-diamond} the two spines have schedules
\[
\begin{aligned}
 A(x,0)&\colon 2x-2,&B(x,0)&\colon 2x-1 &&(1\le x<q),\\
 B(x,-x)&\colon 2x-1,&A(x+1,-x)&\colon 2x &&(1\le x<q/2).
\end{aligned}
\]
Their births follow by induction from the potentials $2x+p$, $-2w-p$
and $2(x+w)+p$, which exclude earlier arrivals and competing parents.
The second list does not repeat the seed.

The upper part $w>0$ maps by $(p,x,w)\mapsto(p,w,x)$ onto $H^+$
truncated by $x'+y'<q$. Its spine schedule is the global schedule
shifted by $-3$ rounds. Its other boundaries are the mirror and A barriers
of Lemma~\ref{lem:isolation}, and returns to its spine are too late.
For every even $m<q$, the inequality $m+\lb(m)\le q$ shows that the
corresponding diamond is retained in full. Coverage gives exactly
\[
 \sum_{\substack{2\le m<q\\m\ \mathrm{even}}}
 z^{2m}\widehat P_{\lb(m)}.
\]

For the lower part apply the inverse cone map
\[
 (p,x,w)\longmapsto
 (P,u,v)=(1-p,-w,x+w+p-1).
\]
The inequalities become $u,v\ge0$ and $2u+v+P<q$.
The seed $A(0,0)$ of this cone corresponds to the external vertex $B(0,0)$,
which precedes the diamond seed by one round. Its next vertex $B(0,0)$
corresponds to the diamond seed $A(1,0)$.
Thus the cone age equals the diamond age plus one. The common horizontal
spine is counted only once, and the two parts are autonomous by the
same spine bounds.

The lower part therefore contributes $(C_q-1)/z$, including the horizontal spine.
The upper part and this cone exhaust the diamond, with disjoint counted vertices. Hence
\[
 \widehat P_q=\sum_{\substack{2\le m<q\\m\ \mathrm{even}}}z^{2m}\widehat P_{\lb(m)}+\frac{C_q-1}{z}.
\]
Induction on $q$ and Definition~\ref{def:P} give $\widehat P_q=P_q$.
The centers with $\lb(m)=q$ are $m=q(2j+1)$ for $j\ge0$, and their seeds occur at
$2q+3+4qj$. Summing the isolated polynomials proves \eqref{eq:F}.
\end{proof}

\subsection{Proof of the generating function identity}

By \eqref{eq:F} and \eqref{eq:G} we have $F=z^2G$, and every $P_q$ has constant term $1$, so the first
nonzero coefficient of $F$ is $F(7)$. Multiply \eqref{eq:filter} by $z^n$ and sum over $n\ge 4$. Since
$F(j)=0$ for $j\le 6$, the four sums over $n\ge 4$ can be extended to all $n$, and
\[
\sum_{n\ge 4}\bigl(2F(n+2)+2F(n+1)+4F(n)+4F(n-1)\bigr)z^n=\bigl(2z^{-2}+2z^{-1}+4+4z\bigr)F=2(1+z+2z^2+2z^3)\,G .
\]
The values $t(1)=1$, $t(2)=2$, and $t(3)=4$ follow from the first three rounds. Therefore
\[
B(z)=z+2z^2+4z^3+\frac{6z^4}{1-z}+2(1+z+2z^2+2z^3)\,G(z).
\]
The rational part equals $z(1+z+2z^2+2z^3)/(1-z)$, and $1+z+2z^2+2z^3=(1+z)(1+2z^2)$. This gives
\eqref{eq:SGprime}, and \eqref{eq:SG} follows from $A(z)=B(z)/(1-z)$. \qed

\section{Reduction of the cones}\label{sec:cones}

By Definition~\ref{def:P} the polynomials $P_q$ are determined by the cones. Most cones can in turn be
reduced to smaller cones.

\subsection{Truncated triangles and the dyadic copy theorem}

For integers $M\ge 1$ and $0\le c\le M$ let
\[
\Omega_{M,c}=\{(p,x,y)\ \mid\ x\ge 0,\ y\ge 0,\ x+y<M,\ 2x+y+p<2M-c\},
\]
and let $L_{M,c}$ be the birth polynomial of the growth on $\Omega_{M,c}$ seeded at $A(0,0)$. We write
$L_M=L_{M,0}$. The domain $\Omega_{M,0}$ is a triangle, and $\Omega_{M,c}$ is obtained by cutting its corner
near the point $x=M$, $y=0$. These domains are closed under predecessors.

\begin{theorem}\label{thm:copy}
Let $m$ be dyadic, let $M$ be a positive multiple of $m$, let $0\le c\le m$ and put $k=M-m$. In the growth on
$\Omega_{M,c}$, the vertices with $x\ge k$ evolve as in the growth on $\Omega_{m,c}$, after the translation
$x\mapsto x-k$ and with all times increased by $2k$. The vertices with $x<k$ evolve as in the growth on
$\Omega_{M,0}$. In particular
\begin{equation}\label{eq:copy}
L_{M,c}-L_{M}=z^{2(M-m)}\bigl(L_{m,c}-L_{m}\bigr).
\end{equation}
\end{theorem}

\begin{proof}
The case $k=0$ is immediate, so let $k>0$. Since $k$ is a multiple of the dyadic integer $m$, every $y$ with
$0\le y<m$ satisfies $k\bitand y=0$ and $(k-1)\bitand y=y$. The vertices of $\Omega_{M,c}$ with $x\ge k$ have
$y<m$, and the translation $x\mapsto x-k$ maps them onto $\Omega_{m,c}$, because the two constraints become
$x+y<m$ and $2x+y+p<2m-c$. Every vertex $A(k,y)$ with $0\le y<m$ belongs to $\Omega_{M,c}$, since
$2k+y\le 2M-m-1<2M-c$. By Lemma~\ref{lem:pascal} it is born at time $2(k+y)$, whatever the value of $c$. The
only edges between the region $x<k$ and the region $x\ge k$ join $B(k-1,y)$ to $A(k,y)$. These B vertices are retained, since $k-1+y<M$ and
$2(k-1)+y+1<2k+y<2M-c$. For $y\ge1$ the
vertex $B(k-1,y)$ is not born at its earliest time $2(k+y)-1$, again by Lemma~\ref{lem:pascal}. Births of a
given type occur at times of a fixed parity, so if $B(k-1,y)$ is born later, its proposal reaches $A(k,y)$ at
time $2(k+y)+2$ or later, when $A(k,y)$ is already born. For $y=0$ the proposal of $B(k-1,0)$ creates
$A(k,0)$ at time $2k$. Hence the region $x\ge k$ receives from outside a single effective proposal, which
makes $A(k,0)$ a seed born at time $2k$. No vertex of this region can be touched before time $2k$. Its growth
is therefore the growth on $\Omega_{m,c}$, translated and delayed by $2k$. In the other direction, every
vertex of $\Omega_{M,0}$ with $x<k$ belongs to $\Omega_{M,c}$, because $2x+y+p\le M+x\le 2M-m-1$. The region
$x<k$ receives proposals from the region $x\ge k$ only through the vertices $A(k,y)$, whose birth times do
not depend on $c$. By induction on time its growth does not depend on $c$ either. Applying the first
statement to $c$ and to $0$ and subtracting gives \eqref{eq:copy}.
\end{proof}

\subsection{Adjacent cones and the Pascal wall}

For $r\ge1$ consider the finite wedge
\[
 \mathcal W_r=\{(p,x,y)\mid x\ge0,\ x+y\ge0,\ 2x+y+p<r\}.
\]
Its seed is $A(0,0)$ and its other A vertices on $x+y=0$ are blocked
initially. All remaining vertices are intact.

\begin{lemma}\label{lem:adjacent}
The birth polynomial of this wedge is
$H_r=C_r+zC_{r-1}-J_{1,r-1}$.
\end{lemma}

\begin{proof}
The potentials $2x+p$ and $2(x+y)+p$ give the lower bound
$\max(2x+p,2(x+y)+p)$ for every first touch.
Along $y=0$, induction therefore gives births $A(x,0)$ at $2x$ and
$B(x,0)$ at $2x+1$ whenever the sites are retained. A return from
$B(x,-1)$ cannot reach $A(x,0)$ before $2x+2$. A return from
$A(x,1)$ cannot reach $B(x,0)$ before $2x+3$.
The regions $y\ge0$ and $y\le0$, $x+y+p\ge1$ are consequently
autonomous inside the wedge, the remaining A boundary being blocked.

The first region is the cone $K_r$. The automorphism
\[
 \Lambda(p,x,y)=(1-p,x+y+p,-y)
\]
maps $K_{r-1}$ onto the second region, because it increases
$2x+y+p$ by one. It sends the seed to $B(0,0)$, born at age $1$.
All other sites of that region are then intact, giving its polynomial
$zC_{r-1}$. Their intersection is the horizontal spine without the seed,
with one birth at each age $1,\ldots,r-1$.
All wedge vertices outside the union are the initially blocked A vertices.
Subtracting the intersection proves the formula.
\end{proof}

\begin{proposition}
Let $M$ be dyadic and $0\le c\le M$. Then
\begin{equation}\label{eq:wall}
 C_{2M-c}=L_{M,c}+z^{2M}H_{M-c}.
\end{equation}
The complete triangle satisfies
\begin{equation}\label{eq:triangle}
 L_M=1+z+2J_{2,2M-1}+2S_{M,2}.
\end{equation}
\end{proposition}

\begin{proof}
In $K_{2M-c}$ the vertices with $x+y<M$ form $\Omega_{M,c}$.
Every retained $A(x,M-x)$ with $x>0$ has $x<M-c$ and receives two
earliest proposals at age $2M$, since $M$ is dyadic. Both parents
belong to the truncated triangle. The vertex is thus permanently blocked.
At the only retained endpoint $A(0,M)$, a possible return reaches
$B(0,M-1)$ after that B has been born at $2M-1$.
The other endpoint $A(M,0)$ is outside $K_{2M-c}$.
These observations prove autonomy of the truncated triangle.

If $c<M$, the region $x+y\ge M$ receives its unique seed
$A(0,M)$ at age $2M$. At this time its retained A vertices on the
wall have just been blocked and all its other vertices are intact,
by the earliest-time bound. Translation $y\mapsto y-M$ gives exactly
$\mathcal W_{M-c}$. The condition $y\ge0$ in the original coordinates is
automatic after this translation, since $x<M-c\le M$ and
the new $y$ is at least $-x$. Lemma~\ref{lem:adjacent} gives its
polynomial $H_{M-c}$ and proves \eqref{eq:wall}. For $c=M$ there is
no residual region and $H_0=0$, giving the same identity.

For $c=0$, the two axes of the triangle meet at $A(0,0)$ and
$B(0,0)$. They contribute $1+z+2J_{2,2M-1}$.
The diagonal reflection separates the remaining births into two equal
halves. Each half contains all the diamonds with even centers
$2\le m<M$ and no others. They fit because $m+\lb(m)\le M$.
Their seeds have cone age $2m+2$. Isolation and coverage therefore
give $2S_{M,2}$, proving \eqref{eq:triangle}.
\end{proof}

\begin{figure}[htbp]
\centering
\includegraphics[width=0.84\linewidth]{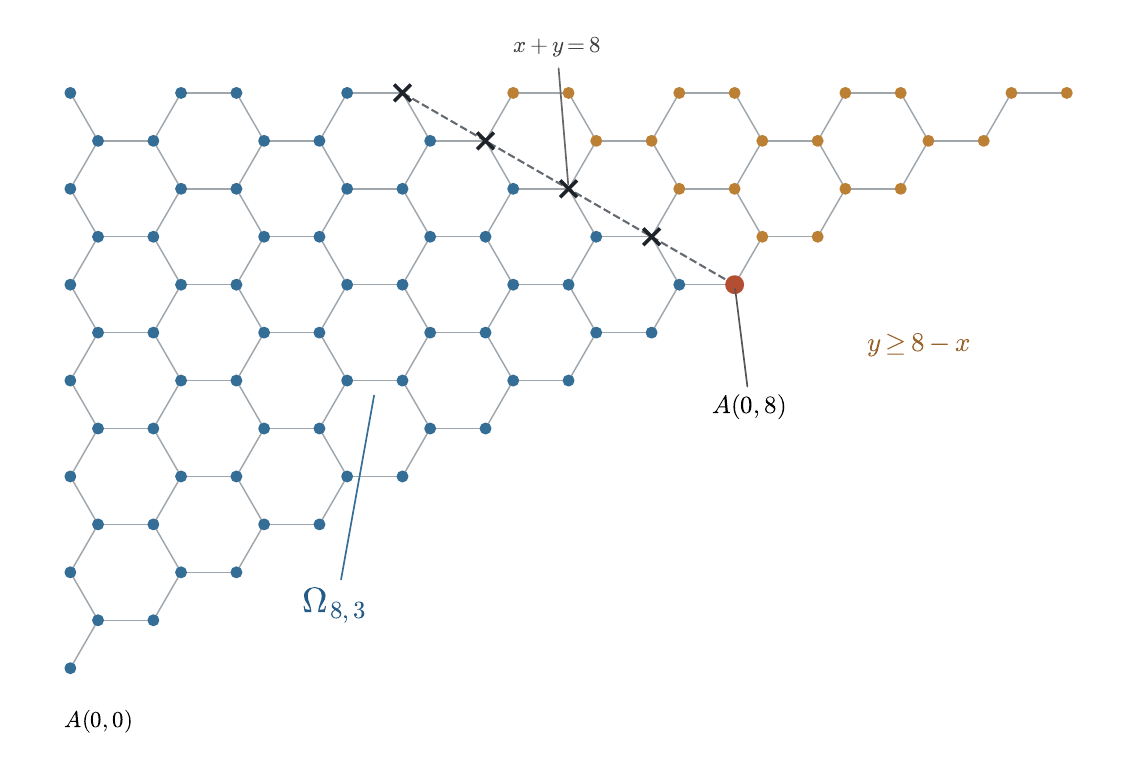}
\caption{The cone $K_{13}$. The blue vertices form $\Omega_{8,3}$. On the line $x+y=8$, the crosses
are A vertices blocked by two arrivals at time $16$. The red vertex $A(0,8)$ is the only seed of the
orange region, whose birth polynomial is $z^{16}H_5$.}
\label{fig:wall}
\end{figure}

Figure~\ref{fig:wall} shows the decomposition when $M=8$ and $c=3$.

\subsection{The reduction}

\begin{proposition}\label{prop:reduction}
Let $r\ge 2$, let $q$ be the smallest dyadic integer with $q\ge r$, and put $c=q-r$ and $M=q/2$. If $c=0$ then
\begin{equation}\label{eq:red0}
C_q=L_M+z^{q}H_M .
\end{equation}
If $c>0$, let $m$ be the smallest dyadic integer with $m\ge c$. Then
\begin{equation}\label{eq:red}
C_r=L_M+z^{q}H_{M-c}+z^{q-2m}\bigl(C_{2m-c}-L_m\bigr)-z^{q}H_{m-c}.
\end{equation}
When $c\le q/4$ every cone on the right-hand side of \eqref{eq:red} has an index smaller than $r$. When
$c>q/4$ one has $m=M$ and \eqref{eq:red} is the trivial identity $C_r=C_r$.
\end{proposition}

\begin{proof}
Since $r>q/2$ we have $0\le c<M$. Identity \eqref{eq:red0} is \eqref{eq:wall} with $c=0$. Let $c>0$. Then
$m\le M$ and both are dyadic, so $M$ is a multiple of $m$ and \eqref{eq:copy} gives
$L_{M,c}=L_M+z^{q-2m}(L_{m,c}-L_m)$. Identity \eqref{eq:wall} at size $m$ gives
$L_{m,c}=C_{2m-c}-z^{2m}H_{m-c}$. Substituting both into \eqref{eq:wall} at size $M$ yields \eqref{eq:red}.
The polynomial $H_{M-c}$ involves the cones of index $M-c$ and $M-c-1$, and $H_{m-c}$ those of index $m-c$
and $m-c-1$, all smaller than $M<r$. Since $m<2c$ we have $2m-c<3c$. If $c\le q/4$ then
$3c\le 3q/4\le q-c=r$, so $2m-c<r$. If $c>q/4$ then $M/2<c<M$, hence $m=M$ and $2m-c=r$, and the terms of
\eqref{eq:red} cancel in pairs.
\end{proof}

With \eqref{eq:triangle} the polynomials $L_M$ and $L_m$ are expressed through diamonds of size smaller than
$M$. A cone is called \emph{direct} when $r\le 3$ or when $c>q/4$. Direct cones are computed by running the
growth of Definition~\ref{def:growth} on $K_r$.

\subsection{Eliminating the diamonds}

Definition~\ref{def:P} also gives a generating function involving the dyadic cones alone.

The product $\mathcal G$ was defined in the introduction. Its coefficient
$2^{s_2(j)}$ counts the odd entries in row $j$ of Pascal's triangle
(OEIS A001316~\cite{OEIS}). Choosing a subset of the nonzero binary digits
of $j$ gives exactly $2^{s_2(j)}$ odd binomial coefficients.

\begin{proposition}
In the ring of formal power series,
\begin{equation}\label{eq:cone-form}
 F(z)=\sum_{q=2,4,8,\ldots}z^{2q+2}\bigl(C_q(z)-1\bigr)\,\mathcal G\bigl(z^{4q}\bigr)
     =\sum_{\substack{m\ge2\\ m\ \mathrm{even}}}2^{s_2(m)-1}z^{2m+2}\bigl(C_{\lb(m)}(z)-1\bigr).
\end{equation}
Moreover, for every dyadic $M\ge1$,
\begin{equation}\label{eq:cone-recursion}
 C_{2M}=1+z+2J_{2,2M-1}
 +\sum_{\substack{2\le m<M\\ m\ \mathrm{even}}}2^{s_2(m)}z^{2m+1}\bigl(C_{\lb(m)}-1\bigr)+z^{2M}H_M .
\end{equation}
\end{proposition}

\begin{proof}
For a dyadic $q$ let $U_q$ be the sum of $z^{2p+3}P_p/(1-z^{4p})$ over the dyadic $p$ with $2\le p\le q$, so
that $U_1=0$. The centers $m<q$ with $\lb(m)=p$ are the odd multiples of $p$ below $q$. Summing the
corresponding geometric progression in \eqref{eq:S} gives $S_{q,0}=z^{-3}(1-z^{2q})U_{q/2}$. With
Definition~\ref{def:P} we obtain
\[
 U_q=\frac{1+2z^{2q}}{1+z^{2q}}\,U_{q/2}+\frac{z^{2q+2}(C_q-1)}{1-z^{4q}} .
\]
By induction, $U_Q$ is the sum over the dyadic $q\le Q$ of $z^{2q+2}(C_q-1)/(1-z^{4q})$ multiplied by the
product of the factors $(1+2z^{2p})/(1+z^{2p})$ over the dyadic $p$ with $q<p\le Q$. Put $w=z^{4q}$. These
factors are $(1+2w^{2^i})/(1+w^{2^i})$ for $0\le i<I$, where $2^I=Q/q$, and the product of the denominators is
$(1-w^{2^I})/(1-w)$. Letting $Q$ tend to infinity gives the first expression in \eqref{eq:cone-form}, since
$F=\lim U_Q$ by \eqref{eq:F}. The second expression follows by writing $m=q(2j+1)$, for which
$s_2(m)=s_2(j)+1$. For a finite $Q=M/2$ the same computation gives
\[
 (1-z^{2M})\,U_{M/2}=\sum_{\substack{2\le m<M\\ m\ \mathrm{even}}}2^{s_2(m)-1}z^{2m+2}\bigl(C_{\lb(m)}-1\bigr),
\]
and the left-hand side is $zS_{M,2}$. Identity \eqref{eq:cone-recursion} is then \eqref{eq:red0} with
\eqref{eq:triangle}.
\end{proof}

Combining \eqref{eq:cone-form} with \eqref{eq:SGprime} proves the second form of Theorem~\ref{thm:series}.
In \eqref{eq:cone-recursion} the polynomial $H_M=C_M+zC_{M-1}-J_{1,M-1}$ involves the cone of size $M-1$,
which is not dyadic. Writing $X_k=C_{2^k}-C_{2^k-1}$ isolates the missing term in this
recursion. Lemma~\ref{lem:adjacent-uniform} bounds its coefficient norm.

\subsection{Layer polynomials}

For a vertex $v=(p,x,y)$ put $g(v)=2x+y+p$. The cone $K_r$ is the set of vertices with $x\ge0$, $y\ge0$ and
$g<r$, and $\Omega_{M,c}$ is the set of vertices of $\Omega_{M,0}$ with $g<2M-c$. The cone $K_{r-1}$ is thus
obtained from $K_r$ by removing the vertices with $g=r-1$, and $\Omega_{M,c+1}$ is obtained from
$\Omega_{M,c}$ by removing the vertices with $g=2M-c-1$. In both cases we call the removed set the
\emph{last layer}. The earliest time of a vertex is $g(v)+y$, so the last layer cannot be reached before
time $r-1$ in the first case and before time $2M-c-1$ in the second. We put
\[
 \partial C_r(z)=z^{-(r-1)}\bigl(C_r(z)-C_{r-1}(z)\bigr)\quad(r\ge1),\qquad \partial C_0=1,
\]
and, for a dyadic integer $M$ and $0\le c<M$,
\begin{equation}\label{eq:kernel}
 E_c(z)=z^{-(2M-c-1)}\bigl(L_{M,c}(z)-L_{M,c+1}(z)\bigr),\qquad e_j(c)=[z^j]E_c(z).
\end{equation}
We call $E_c$ the \emph{kernel} of index $c$. Its coefficients are signed differences of birth counts.
The exponent of $z$ is a delay, counted from the time at which the last layer is reached. By definition
\[
 C_r(z)=\sum_{s=1}^{r}z^{s-1}\,\partial C_s(z),\qquad X_k(z)=z^{2^k-1}\,\partial C_{2^k}(z).
\]

\begin{lemma}\label{lem:kernel-support}
The functions $\partial C_r$ and $E_c$ are polynomials in $z$ with constant term $1$. The kernel $E_c$ does
not depend on the dyadic integer $M>c$. Moreover $\deg\partial C_r\le r(r-1)/2$, and $\deg E_c\le4c^2$ for
$c\ge1$, while $E_0=1$.
\end{lemma}

\begin{proof}
Let $D$ be one of the domains $K_r$ and $\Omega_{M,c}$, let $t_0$ be $r-1$ or $2M-c-1$ accordingly, and let
$D'$ be $D$ without its last layer. A vertex of the last layer has earliest time $t_0+y$. Hence no vertex of the
last layer is born before time $t_0$, and none sends a proposal before time $t_0+1$. The growths on $D$ and on
$D'$ therefore agree on $D'$ up to time $t_0$. The only
vertex of the last layer that can be born at time $t_0$ is the vertex of the axis $y=0$ with $2x+p=t_0$.
It belongs to $D$, and it is born at time $t_0$ by Lemma~\ref{lem:pascal}. The difference of the two birth
polynomials is therefore $z^{t_0}$ plus terms of higher degree.

Let $m$ be the least dyadic integer with $m>c$ and let $M$ be a dyadic multiple of $m$.
Theorem~\ref{thm:copy}, applied to the cuts $c$ and $c+1$, gives
$L_{M,c}-L_{M,c+1}=z^{2(M-m)}(L_{m,c}-L_{m,c+1})$. Thus \eqref{eq:kernel} gives the same polynomial for $M$
and for $m$.

By Remark~\ref{rem:finite-domain}, the $r(r+1)/2$ vertices of $K_r$ give the bound on
$\deg\partial C_r$. The triangle $\Omega_{m,0}$ has $m(m+1)$ vertices, so $\deg E_c\le m^2-m+c$. For $c\ge1$ we
have $m\le2c$. For $c=0$ we have $m=1$, $L_{1,0}=1+z$ and $L_{1,1}=1$.
\end{proof}

For example
\[
 E_0=1,\qquad E_1=E_2=1+z,\qquad E_3=1+z^2+2z^3,\qquad E_4=1+z^2+z^3-z^5,
\]
and $\partial C_1=1$, $\partial C_2=1+z$, $\partial C_3=1+z+z^2$, $\partial C_4=1+z^2+z^3$.
For $n\ge0$ let $J_n$ be the largest degree of the kernels $E_0,\ldots,E_n$. By the lemma, $J_0=0$ and
$J_n\le4n^2$.

\subsection{Differences of consecutive cones}

\begin{proposition}\label{prop:transfer}
Let $M$ be dyadic and let $M<r\le2M$. Then
\begin{equation}\label{eq:transfer}
 \partial C_r=E_{2M-r}+z^{M}\bigl(\partial C_{r-M}+\partial C_{r-M-1}-1\bigr).
\end{equation}
\end{proposition}

\begin{proof}
Put $c=2M-r$, so that $0\le c<M$. Identity \eqref{eq:wall} for the cuts $c$ and $c+1$ gives
\[
 C_r-C_{r-1}=L_{M,c}-L_{M,c+1}+z^{2M}\bigl(H_{M-c}-H_{M-c-1}\bigr).
\]
Let $\rho=M-c=r-M$. For $\rho\ge2$ the definition of $H_\rho$ gives
$H_\rho-H_{\rho-1}=(C_\rho-C_{\rho-1})+z(C_{\rho-1}-C_{\rho-2})-z^{\rho-1}$, which is
$z^{\rho-1}(\partial C_\rho+\partial C_{\rho-1}-1)$. For $\rho=1$ we have $H_1-H_0=1$, which is the same
expression because $\partial C_1=\partial C_0=1$. The first difference is $z^{r-1}E_c$ by \eqref{eq:kernel}.
Dividing by $z^{r-1}$ gives \eqref{eq:transfer}, since $2M+\rho-1-(r-1)=M$.
\end{proof}

For instance $\partial C_8=E_0+z^4(\partial C_4+\partial C_3-1)=1+z^4+z^5+2z^6+z^7$. At $z=1$, this identity relates the total numbers of births.
Its coefficient norm controls the error made by replacing a cone by its neighbor.
Its successive substitutions give the decomposition of Section~\ref{sec:windows}.

\begin{lemma}\label{lem:adjacent-uniform}
For every integer $k\ge1$,
\[
 \|C_{2^k}-C_{2^k-1}\|=O(k^2 2^k),
\]
with an absolute implied constant. The same bound holds for the absolute value of any partial sum of these coefficients.
\end{lemma}

\begin{proof}
The kernel $E_c$ is computed in a triangle of size at most $2\max(c,1)$, so $\|E_c\|=O((c+1)^2)$.
Apply \eqref{eq:transfer} repeatedly to $\partial C_{2^k}$. At a node with cone index $2^a-c$, continue
only while $c+1\le2^{a-2}$. The depth is at most $k$, and every cut index is at most $k$.
Each internal node contributes $O((k+1)^2)$ to the norm. At a leaf, the cone index is smaller than
$4(k+1)$, so the finite-domain bound gives the same estimate. There are at most $2^k$ leaves and fewer
than $2^k$ internal nodes. Multiplication by a monomial preserves the norm, which proves the result.
Explicit constants are given in Appendix~\ref{app:bounds}.
\end{proof}

\section{Coefficient extraction and complexity}

Throughout this section $n\ge 4$ is the queried index and $\ell=\lfloor\log_2(n+3)\rfloor+1$. By
\eqref{eq:filter} and by summation, with $\Sigma F(j)=\sum_{i\le j}F(i)$,
\begin{equation}\label{eq:an}
a(n)=6n-11+2\,\Sigma F(n+2)+2\,\Sigma F(n+1)+4\,\Sigma F(n)+4\,\Sigma F(n-1).
\end{equation}
It is therefore enough to compute a coefficient $F(t)$ or a partial sum $\Sigma F(t)$ for $t\le n+2$. The
value $t(n)$ is obtained from \eqref{eq:filter} with four coefficients. Only
the diamonds $P_q$ with $2q+3\le n+2$ are involved, and their sizes are $q=2^i$ with $i<\ell$.

\subsection{A graph of expressions}

The algorithm builds a finite directed acyclic graph. Its nodes are the polynomials $P_q$ and $L_M$ for
dyadic indices, and the polynomials $C_r$ and $H_r$ that are reached from them. Each node is written
as a signed sum of \emph{blocks}
\begin{equation}\label{eq:block}
R(z)=\sum_{j=0}^{N-1}z^{\alpha+\delta j}\,Q(z),
\end{equation}
where $Q$ is another node, an interval $J_{u,v}$ or the constant $1$, and where $\alpha$, $\delta$, and $N$
are nonnegative integers. For $q\ge2$ dyadic, \eqref{eq:diamond-cone}, \eqref{eq:red0}, and
\eqref{eq:triangle} give the rule
\begin{equation}\label{eq:diamond-rule}
 P_q=1+2J_{1,q-2}+S_{q,0}+2S_{q/2,1}+z^{q-1}H_{q/2}.
\end{equation}
All its shifts are nonnegative, including when $q=2$. The other rules are the definition of $H_r$,
\eqref{eq:triangle} for $L_M$, and \eqref{eq:red0} and \eqref{eq:red} for cones that are not direct. The direct cones are leaves. Their coefficients are computed by running the growth and are
stored as lists. Every node is created once and shared by all the expressions that refer to it.

The sums $S_{B,d}$ have about $B/2$ terms. To group them into blocks, partition the even integers
$m$ according to $p=\lb(m)$. The integers $m<B$ with $\lb(m)=p$ are the numbers $p(2j+1)$ with
$0\le j<N_{B,p}$, where $N_{B,p}=1+\lfloor (B-1-p)/(2p)\rfloor$. Hence
\begin{equation}\label{eq:compress}
S_{B,d}=\sum_{\substack{p=2,4,8,\ldots\\ p<B}}\ \sum_{j=0}^{N_{B,p}-1}z^{2p+d+4pj}\,P_p ,
\end{equation}
which is a sum of fewer than $\log_2 B$ blocks. The series $F$ is treated in the same way. For a query at
index $t$ only the first $1+\lfloor (t-2q-3)/(4q)\rfloor$ copies of $P_q$ in \eqref{eq:F} can contribute.

To each node $X$ we attach its total mass $\mu(X)=X(1)$ and an integer $b(X)$ that bounds its degree. Both are
computed from the rules. The bound $b$ of a sum is the largest of the bounds of its shifted terms, without
regard to cancellations. Zero terms are omitted, and a direct cone is assigned its degree.

\subsection{The cones that occur}

For $r\ge 1$ let $\lambda(r)$ be the least integer with $2^{\lambda(r)}\ge r$, let $c(r)=2^{\lambda(r)}-r$, and call
$\rho(r)=\lambda(r)+c(r)$ the \emph{rank} of $r$. If $r=2^b-e\ge 1$ with integers $b\ge 0$ and $e\ge 0$, then
$\rho(r)\le b+e$. Indeed, either $e<2^{b-1}$ and $b=\lambda(r)$, or $r\le 2^{b-1}$ and $r=2^{b-1}-(e-2^{b-1})$
with a smaller exponent and a smaller deficit, and one repeats the argument.

\begin{lemma}\label{lem:rank}
Every cone $C_r$ of the graph built for the diamonds $P_{2^i}$ with $i\le s$ satisfies $\rho(r)\le s$.
Consequently the graph contains $O(s^2)$ cones, and every direct cone has an index $r<3s$ or $r\le 3$.
\end{lemma}

\begin{proof}
The diamond $P_q$ with $q=2^i$ refers to $H_{q/2}$, hence to the cones of index $2^{i-1}$ and $2^{i-1}-1$,
whose ranks are at most $i$. It also refers to diamonds of smaller size. The triangle $L_M$ refers only to
diamonds of size smaller than $M$. Consider now a cone $C_r$ that is not direct, with $a=\lambda(r)$, $c=c(r)$,
$q=2^a$ and $M=q/2$, so that $c\le q/4$. Through $H_{M-c}$ it refers to the cones of index $2^{a-1}-c$ and
$2^{a-1}-(c+1)$, whose ranks are at most $a+c$. When $c>0$, let $m=2^{b}$ be the smallest dyadic integer with
$m\ge c$. Then $b\le a-2$. The cone $C_r$ refers to the cones of index $2^{b+1}-c$, $2^{b}-c$ and
$2^{b}-(c+1)$, whose ranks are at most $b+1+c\le a+c$. It also refers to $L_M$ and $L_m$, hence to diamonds
of size at most $2^{a-2}$. Along the references between cones the rank never increases, and the cones
referred to by a diamond of size $2^i$ have rank at most $i$. This proves the first claim. A cone of rank at
most $s$ is determined by a pair $(a,c)$ with $a+c\le s$, and there are $O(s^2)$ such pairs. If a cone with
$r\ge 4$ is direct then $c>2^{a}/4$, so $r=2^a-c<3c\le 3(s-a)<3s$.
\end{proof}

\subsection{Extraction of a coefficient and of a partial sum}

For a node $Q$ write $Q_t$ for the coefficient of $z^t$ and $Q_{\le t}=\sum_{i\le t}Q_i$. These numbers vanish
for $t<0$, and for $t\ge b(Q)$ one has $Q_{\le t}=\mu(Q)$. Consider a block $R$ as in \eqref{eq:block}, with
$\delta\ge 1$, and let $b=b(Q)$. For a block with a single copy we take $\delta=1$. For $t\ge\alpha$ put
\begin{gather*}
j_0=\max\Bigl(0,\Bigl\lceil\frac{t-\alpha-b}{\delta}\Bigr\rceil\Bigr),\qquad
j_1=\min\Bigl(N-1,\Bigl\lfloor\frac{t-\alpha}{\delta}\Bigr\rfloor\Bigr),\\
f=\max\Bigl(0,\min\Bigl(j_1+1,\,1+\Bigl\lfloor\frac{t-\alpha-b}{\delta}\Bigr\rfloor\Bigr)\Bigr).
\end{gather*}
The copies with $j<j_0$ end before $t$ and the copies with $j>j_1$ start after $t$. Therefore
\begin{equation}\label{eq:extract}
R_t=\sum_{j=j_0}^{j_1}Q_{t-\alpha-\delta j},\qquad
R_{\le t}=f\,\mu(Q)+\sum_{j=f}^{j_1}Q_{\le t-\alpha-\delta j}.
\end{equation}
In the second formula the copies with $j<f$ lie entirely below $t$ and are replaced by their mass. For
$t<\alpha$ both quantities vanish. A query at a node is answered by applying \eqref{eq:extract} to each of its
blocks, with the appropriate signs. The answers are stored in a table indexed by the node, the index, and the
type of query, so that no query is evaluated twice. The number of terms in each sum of \eqref{eq:extract} is
at most $1+b/\delta$.

\subsection{Supports and residues}

The \emph{scale} $\sigma(X)$ of a node is $q$ for $P_q$ and $L_q$, and it is the smallest dyadic integer
$q\ge r$ for $C_r$ and $H_r$. In every rule, the nodes that appear on the right-hand side have a scale
at most equal to the scale of the node that is defined.

For $k\ge1$, let $\tau_k$ be the maximum of $\deg C_r-2r+2$ over direct cones of rank at most $k$.
The axis vertex $A(0,r-1)$ is born at time $2r-2$, so $\tau_k\ge0$.
Remark~\ref{rem:finite-domain} and Lemma~\ref{lem:rank} give $\tau_k\le5k^2$.

\begin{lemma}\label{lem:degrees}
Let $k\ge1$. For every cone index $r$ with $\rho(r)\le k$, and every dyadic $q=2^j$ with $1\le j\le k$,
\[
 b(C_r)\le2r-2+\tau_k,\qquad b(P_q)\le2q-3+\tau_k,\qquad b(L_q)\le2q-1+\tau_k.
\]
Consequently every node $X$ in the graph for the queried index $n$ satisfies
$\deg X\le b(X)\le2\sigma(X)+5\ell^2$.
\end{lemma}

\begin{proof}
Put $\tau=\tau_k$ and use induction along the rules, starting from the direct cones, whose bounds are
their degrees. If the bounds hold for the cones in $H_r$, its defining terms give
$b(H_r)\le2r-2+\tau$. A term $z^{2m+2}P_p$ of $L_M$ has $p=\lb(m)$ and $m+p\le M$.
Together with the axes, this gives $b(L_M)\le2M-1+\tau$.

For a cone that is not direct, use the notation of Proposition~\ref{prop:reduction}.
All the smaller cones have rank at most $k$, and the diamonds in $L_M$ and $L_m$ have size at most $q/4$.
The five terms of \eqref{eq:red} therefore have bounds at most
\[
 \max\bigl(4M-2c-2,\ 2M+2m-2c-2,\ 2M-1\bigr)+\tau
 \le2r-2+\tau,
\]
because $m\le M$ and $2c+1\le2M$. The case \eqref{eq:red0} follows from the same bounds with $c=0$.
In \eqref{eq:diamond-rule}, a term of $S_{q,0}$ has bound at most $2(m+p)-3+\tau\le2q-3+\tau$.
The terms of $S_{q/2,1}$ have bound at most $q-2+\tau$, and the last term has bound
$q-1+q-2+\tau$. The constant and the interval satisfy the same bound. This proves the asserted
inequalities without using cancellation. Lemma~\ref{lem:rank} and $\tau_\ell\le5\ell^2$ give the conclusion.
\end{proof}

The horizontal spine in the diamond also gives $\deg P_q\ge2q-3$.

\begin{lemma}\label{lem:residue}
Fix the queried index $n$. There is a constant $K_0$ such that every index $t$ queried at a node of scale
$\sigma$ during the evaluation of \eqref{eq:an} satisfies $t\equiv n+\varepsilon\pmod{\sigma}$ for some
integer $\varepsilon$ with $|\varepsilon|\le K_0\ell$.
\end{lemma}

\begin{proof}
A query at a node is generated by a chain of blocks that starts at the root $F$. Passing through a block
\eqref{eq:block} replaces the index $t$ by $t-\alpha-\delta j$. We inspect the shifts $\alpha+\delta j$ of
the rules. In \eqref{eq:compress} the shift is $2p+d+4pj$ with $d\in\{0,1,2\}$, and the child $P_p$ has scale
$p$. In \eqref{eq:diamond-rule} the last term has shift $q-1$ and its child $H_{q/2}$ has scale $q/2$. In
\eqref{eq:red0} and \eqref{eq:red} the shifts are $0$, $q$, and $q-2m$, the children with shift $q$ have a
scale that divides $q$, and the children with shift $q-2m$ have a scale that divides $2m$, which divides
$q-2m$. In the definition of $H_r$ the shifts are $0$ and $1$. In \eqref{eq:F} the shift is $2q+3+4qj$ and
the child has scale $q$. In each case the shift is a multiple of the scale of the child plus an integer of
absolute value at most $3$. Scales are powers of two and do not increase along a chain, so the scale of the
current node divides the scale of every earlier child. Modulo the scale of the current node, the index $t$
therefore differs from the root index by a sum of at most three units per block of the chain. The root
indices are $n+2$, $n+1$, $n$, and $n-1$. The length of a chain is also bounded. Along a chain the scale
is at least halved at each block, except on a passage from $H_r$ to a cone. Two such passages cannot be consecutive. Since there are at most $\ell$ scales, a chain has $O(\ell)$ blocks.
\end{proof}

\subsection{Proof of the direct access theorem}

The terms $a(0),\ldots,a(3)$ are returned directly, so let $n\ge 4$. By Lemma~\ref{lem:rank} the graph has
$O(\ell^2)$ cones, hence $O(\ell^2)$ nodes, since only a bounded number of nodes $H_r$ occur for each cone and
there are $O(\ell)$ diamonds and triangles. By \eqref{eq:compress} each node is a sum of $O(\ell)$ blocks. The
direct cones have $O(\ell)$ distinct indices smaller than $3\ell$ and $O(\ell^2)$ vertices each. Every vertex
receives at most three proposals and emits at most once, so all the direct cones are computed with $O(\ell^3)$
elementary operations on vertices.

We count the entries of the table of answers. At a node of scale $\sigma$ the queried indices lie in the
interval from $0$ to $b$, with $b\le2\sigma+5\ell^2$ by Lemma~\ref{lem:degrees}, because the other indices
are answered by $0$ or by the mass. If $\sigma\le \ell$ this interval contains $O(\ell^2)$ integers. If $\sigma>\ell$,
Lemma~\ref{lem:residue} leaves $O(\ell)$ residues modulo $\sigma$, and each residue has $O(1+\ell^2/\sigma)$
representatives in the interval, which gives $O(\ell^2)$ indices again. With two types of queries the table has
$O(\ell^4)$ entries.

An entry is computed by applying \eqref{eq:extract} to $O(\ell)$ blocks. For a block with child $P_p$ the step
is $\delta=4p$ and the number of terms is at most $1+b(P_p)/(4p)=O(\ell^2)$ by Lemma~\ref{lem:degrees}. The other
blocks have a single copy. An entry therefore costs $O(\ell^3)$ arithmetic operations and look-ups, and the
whole evaluation costs $O(\ell^7)$ of them. With a balanced dictionary each look-up costs $O(\log \ell)$
comparisons, giving $O(\ell^8)$ operations in total.

All the integers that occur are indices, numbers of copies, or signed sums of numbers of births in domains
with $O(n^2)$ vertices. They have $O(\ell)$ bits. The algorithm is therefore polynomial in $\ell$ in bit
complexity as well. With elementary integer arithmetic, $O(\ell^{10})$
is a sufficient bound in a model with addressed memory. Simulation on a
Turing machine preserves polynomial time, without asserting the same exponent.
Consequently the maps $n\mapsto a(n)$ and $n\mapsto t(n)$ belong to FP. \qed

\begin{remark}
The proof uses neither the identity $\deg P_q=2q-3$ nor any bound on the time at which a cone becomes
extinct other than its number of vertices.
\end{remark}

\section{Computations}\label{sec:exp}

The supplementary code contains two separately written implementations,
both using integer arithmetic and neither reading a table of birth
polynomials. The first builds the expanded cone identities. The second
uses the shorter triangle reduction \eqref{eq:red}, an event queue for
the direct cones, and the factorized generating function
\eqref{eq:SGprime} for the global count. They share the mathematical
identities, but not their implementation.

Both programs agree with direct simulation through $n=1024$, including the OEIS values \cite{OEIS}.
They also pass coefficient and partial-sum checks against the reference diamonds through $P_{65536}$.
The generating function agrees with the global filter through $n=262144$. The README and recorded
outputs give the query counts and the reproduction commands.

The two programs agree on $a(n)$ and $t(n)$ at each of the six large indices in Table~\ref{tab:times}.
For example,
\[
 a(131072)=17\,727\,196\,305,\qquad a(262144)=70\,909\,299\,829,
\]
and
\[
 a(10^{12})=897\,974\,080\,671\,417\,343\,896\,245,\qquad
 t(10^{12})=148\,595\,720\,694.
\]
The value of $a(10^{100})$ has 200 decimal digits. The supplementary
JSON files contain the full integer outputs and the check results.

\begin{table}[htbp]
\centering
\begin{tabular}{rrrrr}
\toprule
Index $n$ & Program I (s) & Program II (s) & Largest cone & I nodes\\
\midrule
$262\,144$ & $0.005$ & $0.004$ & $23$ & $277$\\
$10^{12}$ & $0.060$ & $0.041$ & $47$ & $1\,754$\\
$10^{20}$ & $0.516$ & $0.265$ & $95$ & $5\,454$\\
$10^{30}$ & $1.687$ & $1.105$ & $191$ & $12\,864$\\
$10^{50}$ & $9.081$ & $6.807$ & $383$ & $37\,835$\\
$10^{100}$ & $102.503$ & $63.932$ & $767$ & $157\,240$\\
\bottomrule
\end{tabular}
\caption{Observed wall-clock times for single-process Python implementations.
These are individual runs, not estimates of the optimal complexity exponent.}
\label{tab:times}
\end{table}

The large-index comparisons check the implementations against each other.
They do not constitute independent simulations of the full geometric
structure at those indices, nor are they used in the proofs.

\section{Asymptotic structure}

The normalized masses of the diamonds converge to a constant $\kappa$. The cumulative birth counts,
rescaled on dyadic intervals, converge to a continuous profile. Two functional equations determine
this profile and its extreme values.

\subsection{A computable density for the diamonds}

Write $\mu_q=P_q(1)$ for dyadic $q\ge2$.

\begin{theorem}\label{thm:density}
There is a positive real number $\kappa$ such that, as $k\to\infty$,
\[
 \mu_{2^k}=\kappa4^k+O(k^2 2^k).
\]
The constant is computable, with decimal expansion
\[
 \kappa=0.687915685568298445216698239913\ldots.
\]
\end{theorem}

\begin{proof}
Write $c_r=C_r(1)$ and $d_k=\partial C_{2^k}(1)=c_{2^k}-c_{2^k-1}$.
Lemma~\ref{lem:adjacent-uniform} gives $d_k=O(k^2 2^k)$.
For dyadic $q\ge2$, put $w_q=\sum_{p=2,4,\ldots,q/2}\mu_p/p$.
Evaluating \eqref{eq:diamond-cone}, \eqref{eq:triangle}, and \eqref{eq:red0} at $z=1$ gives
\[
 q w_q=2\mu_q-2c_q+2,\qquad
 c_{2q}=3q+1+2\mu_q-(c_q-c_{q-1}).
\]
The initial case $q=2$ follows from $P_2=1+z$ and $C_2=1+z+z^2$.
With $\delta_k=\mu_{2^{k+1}}-4\mu_{2^k}$, these identities imply
\begin{equation}\label{eq:mass-increments}
 \delta_1=2,\qquad \delta_k=\delta_{k-1}+2d_{k-1}-d_k\quad(k\ge2).
\end{equation}
Telescoping gives $\delta_k=O(k^2 2^k)$. Therefore
\[
 \kappa=\frac12+\sum_{k\ge1}\frac{\delta_k}{4^{k+1}}
\]
converges absolutely, and its tail gives the stated estimate.
Appendix~\ref{app:bounds} gives an explicit error bound, proves positivity from a computed value,
and certifies the displayed decimal digits.
\end{proof}

\subsection{A rapidly convergent series for the density}

The transfer identity also gives a series whose first $N$ terms use cones of size $O(N)$.
For $j\ge0$, put
\[
 Q_j=2^{\lceil\log_2(j+2)\rceil},\qquad
 \epsilon_j=\mathbf1_{j=Q_j/2-1},\qquad
 \alpha_j=\frac{3\partial C_{Q_j-j}(1)-1+E_j(1)+9\epsilon_j}{Q_j^2}.
\]
Thus the row sums $E_j(1)$ of the kernels determine the density, while their columns determine the
windows in Section~\ref{sec:windows}.

\begin{theorem}\label{thm:ternary}
The density has the absolutely convergent representation
\[
 \kappa=1-\frac2{27}\sum_{j\ge0}\frac{\alpha_j}{3^j},
 \qquad
 \left|\kappa-1+\frac2{27}\sum_{j=0}^{N-1}\frac{\alpha_j}{3^j}\right|
 \le\frac4{3^N}\quad(N\ge1).
\]
The first $N$ terms use only cones of index at most $4N$.
\end{theorem}

\begin{proof}
Let $d_k=c_{2^k}-c_{2^k-1}$ and
$\delta_k=\mu_{2^{k+1}}-4\mu_{2^k}$. The recurrence \eqref{eq:mass-increments}, with $d_1=\delta_1=2$, gives
\[
 (1-t)\Delta(t)=4t+(2t-1)D(t),\qquad
 D(t)=\sum_{k\ge1}d_kt^k,\quad
 \Delta(t)=\sum_{k\ge1}\delta_kt^k.
\]
The estimates already proved justify evaluation at $t=1/4$. Since
$\kappa=1/2+\Delta(1/4)/4$, we obtain
\begin{equation}\label{eq:kappa-dyadic-series}
 \kappa=\frac{5-D(1/4)}6
 =\frac56-\frac16\sum_{k\ge1}\frac{c_{2^k}-c_{2^k-1}}{4^k}.
\end{equation}

For $a,j\ge0$ let $Y(a,j)=c_{2^a-j}-c_{2^a-j-1}$ if $j<2^a$, and let
$Y(a,j)=0$ otherwise. Evaluating \eqref{eq:transfer} at $z=1$ gives
\begin{equation}\label{eq:ternary-transfer}
 Y(a,j)=Y(a-1,j)+Y(a-1,j+1)+g_j,\qquad g_j=E_j(1)-1,
\end{equation}
whenever $2^{a-1}\ge j+2$.

For convergence, expand only while $j+1\le2^{a-2}$. In the tree starting
at $(a,j)$ every corner index is at most $j+a$. The estimates in Appendix~\ref{app:bounds} give
\[
 |Y(a,j)|\le40(a+j+1)^2 2^a,\qquad
 U_j:=\sum_{a\ge0}Y(a,j)4^{-a},\qquad
 |U_j|\le80(j^2+4j+6).
\]
Thus these series converge absolutely.

Write $Q_j=2^h$ and $w_j=c_{Q_j-j}-c_{Q_j-j-1}$. Below level $h$ the
only possible nonzero value is $Y(h-1,j)=1$, and it occurs exactly when
$\epsilon_j=1$. At level $h$ the value is $w_j$, and at every level $a\ge h+1$
recurrence~\eqref{eq:ternary-transfer} applies. Therefore, with
$Y_j(t)=\sum_{a\ge0}Y(a,j)t^a$,
\[
 (1-t)Y_j(t)-tY_{j+1}(t)
 =w_jt^h+\epsilon_j(1-t)t^{h-1}+\frac{g_jt^{h+1}}{1-t}.
\]
At $t=1/4$ this gives
\[
 U_j=\frac13U_{j+1}+\frac49\alpha_j.
\]
Since $U_j=O(j^2)$, iteration yields
$U_0=(4/9)\sum_{j\ge0}\alpha_j3^{-j}$. Since $Y(0,0)=1$, we have
$D(1/4)=U_0-1$. Substitution in \eqref{eq:kappa-dyadic-series} proves the formula.

The bounds on $E_j$ in Appendix~\ref{app:bounds} and $0\le c_r\le r(r+1)/2$ imply
\[
 |g_j|\le29(j+1)^2\le29Q_j^2,\qquad
 \frac{|w_j|}{Q_j^2}\le\frac34,\qquad
 |\alpha_j|\le29+\frac94+\frac94<36.
\]
The tail is at most $(2/27)\cdot36\sum_{j\ge N}3^{-j}=4/3^N$.
For $j<N$, the kernel $E_j$ is computed with a dyadic size $m\le2\max(j,1)$.
Identity \eqref{eq:wall} expresses its two triangle masses through cones of index at most $2m$.
Also $Q_j<2(j+2)\le4N$, so every required cone index is at most $4N$.
\end{proof}

The supplied program certifies $105$ decimal digits from $224$ terms, using cones of
index at most $384$, with every cone mass checked by two implementations.

\subsection{Block identities}

Identity \eqref{eq:diamond-cone} relates the global count to a dyadic cone.

For $r\ge0$ and an integer $t$, let $c_r(t)$ be the coefficient of $z^t$ in
$C_r$ and let $\Gamma_r(t)=\sum_{i\le t}c_r(i)$. Let $\mathcal H(n)=\sum_{i\le n}F(i)$. These quantities vanish
for negative arguments. By \eqref{eq:F}, the coefficient $F(n)$ is the sum of the coefficients of
$z^{n-2m-3}$ in $P_{\lb(m)}$ over the even integers $m\ge2$. The integer $m$ is the center of a diamond whose
seed appears at round $2m+3$. For a dyadic integer $B$ let $F^{<B}(n)$ be the same sum restricted to $m<B$, and
let $\mathcal H^{<B}(n)=\sum_{i\le n}F^{<B}(i)$. Then $F^{<B}(n)=F(n)$ for $n\le2B+2$. The function
$\mathcal H^{<B}$ is nondecreasing and its limit is
\[
 \mathcal M_B=\sum_{\substack{2\le m<B\\ m\ \mathrm{even}}}\mu_{\lb(m)}
    =\sum_{p=2,4,\ldots,B/2}\frac{B}{2p}\,\mu_p .
\]

\begin{lemma}\label{lem:block}
Let $N=2q$ with $q\ge2$ dyadic. For $3\le j\le N$,
\begin{equation}\label{eq:block-rule}
 F(N+j)=2F(j)+c_q(j-2)+F^{<q}(N+j),
\end{equation}
and $F(N+j)=F^{<q}(N+j)$ for $j=1$ and $j=2$. Consequently, for $0\le J\le N$,
\begin{equation}\label{eq:block-sum}
 \mathcal H(N+J)=\mathcal H^{<q}(N+J)+2\mathcal H(J)+\max\bigl(\Gamma_q(J-2)-1,\,0\bigr).
\end{equation}
\end{lemma}

\begin{proof}
Let $1\le j\le N$. In the sum that defines $F(N+j)$, a center $m\ge2q$ has $2m+3>N+j$ and does not
contribute. The centers $m<q$ contribute $F^{<q}(N+j)$. A center with $q<m<2q$ is $m=q+m'$ with $0<m'<q$. Then
$\lb(m)=\lb(m')$ and the exponent is $j-2m'-3$. These centers therefore contribute $F^{<q}(j)$, which equals
$F(j)$ because $j\le2q+2$. The center $m=q$ contributes the coefficient of $z^{j-3}$ in $P_q$. By
Definition~\ref{def:P} this coefficient is zero for $j\le2$. For $j\ge3$ it is the sum of $c_q(j-2)$ and
of the coefficient of $z^{j-3}$ in $S_{q,0}$. By \eqref{eq:S} the latter is $F^{<q}(j)=F(j)$. This proves
\eqref{eq:block-rule}. Since $F(1)=F(2)=0$ and $\mathcal H(N)=\mathcal H^{<q}(N)$, summation over
$1\le j\le J$ gives \eqref{eq:block-sum}. Its last term is $\sum_{j=3}^{J}c_q(j-2)=\Gamma_q(J-2)-c_q(0)$ when
$J\ge2$.
\end{proof}

\begin{lemma}\label{lem:cone-coeff}
Let $M\ge1$ be dyadic and $N=2M$. For every $t\ge0$,
\begin{equation}\label{eq:cone-coeff}
 c_N(t)=\varepsilon_N(t)+2F^{<M}(t+1)+c_M(t-N)+c_{M-1}(t-N-1)-\mathbf 1_{1\le t-N\le M-1},
\end{equation}
where $\varepsilon_N(t)$ equals $1$ for $t\le1$, equals $2$ for $2\le t\le N-1$ and equals $0$ for $t\ge N$.
\end{lemma}

\begin{proof}
This is \eqref{eq:red0} with \eqref{eq:triangle} and the definition of $H_M$. By \eqref{eq:S} the coefficient
of $z^t$ in $S_{M,2}$ is $F^{<M}(t+1)$.
\end{proof}

\subsection{The limit profile}

The block identity and the reduction of a dyadic cone have continuous counterparts. Let $\mathcal X$ be the set
of pairs $(f,g)$ of continuous functions on $[1,2]$ with
\[
 f(1)=\tfrac18,\qquad f(2)=\tfrac12,\qquad g(1)=\tfrac14,\qquad g(2)=\tfrac12 .
\]
For such a pair, extend $f$ to $[0,\infty)$ by $f(2u)=4f(u)$ and $f(0)=0$, and extend $g$ to $[0,1]$ by
$g(y)=2f(y)$. Both extensions are continuous. Define $\mathcal T(f,g)=(f^*,g^*)$ by
\begin{equation}\label{eq:profile-operator}
 f^*(1+x)=\tfrac18+2f(x)+\tfrac14\,g(2x),\qquad
 g^*(1+x)=\tfrac14+\tfrac12\,g(2x)\qquad(0\le x\le1).
\end{equation}
The four endpoint values are preserved, so $\mathcal T$ maps $\mathcal X$ to itself.

\begin{lemma}
The operator $\mathcal T$ has a unique fixed point $(\mathcal S,\mathcal C)$ in $\mathcal X$. The functions $\mathcal S$ and $\mathcal C$ are
nondecreasing. With $\mathcal S(u)=\mathcal C(u)=0$ for $u\le0$, they satisfy
\begin{equation}\label{eq:profile-single}
 \mathcal C(y)=\mathcal S(2+y)-2\mathcal S(y)-\tfrac12\quad(0\le y\le2),\qquad
 \mathcal S(u)=\tfrac18+\mathcal S(2u-2)-\mathcal S(2u-3)\quad(1\le u\le2).
\end{equation}
\end{lemma}

\begin{proof}
Let $(f_1,g_1)$ and $(f_2,g_2)$ be in $\mathcal X$ and let $\delta$ be the larger of the uniform distances
between $f_1$ and $f_2$ and between $g_1$ and $g_2$ on $[1,2]$. By scaling, the extensions of $f_1$ and $f_2$
differ by at most $\delta/4$ on $[0,1]$, and the extensions of $g_1$ and $g_2$ differ by at most $\delta$ on
$[0,2]$. By \eqref{eq:profile-operator} the images differ by at most $3\delta/4$. Thus $\mathcal T$ is a
contraction for the uniform distance, and it has a unique fixed point $(\mathcal S,\mathcal C)$. The pairs of nondecreasing
functions form a closed subset of $\mathcal X$ which is preserved by $\mathcal T$, because the extensions of
nondecreasing functions are nondecreasing. This subset contains the pair of affine functions with the
prescribed endpoint values, so it contains the fixed point.

The first equation of \eqref{eq:profile-operator} at $x=y/2$, multiplied by~$4$, is the first identity of
\eqref{eq:profile-single}. Let $0\le w\le1$. Substituting this identity at $1+w$ and at $2w$ in
$\mathcal C(1+w)=\frac14+\frac12\mathcal C(2w)$ and using the scaling of $\mathcal S$ gives $\mathcal S(3+w)=4\mathcal S(1+w)-4\mathcal S(w)+\frac12$. With
$u=(3+w)/2$ this is the second identity of \eqref{eq:profile-single} for $3/2\le u\le2$. For $1\le u\le3/2$
the second identity reads $\mathcal S(u)=\frac18+4\mathcal S(u-1)$, which is the first equation of \eqref{eq:profile-operator}
with $g(2x)=2\mathcal S(2x)=8\mathcal S(x)$.
\end{proof}

The second identity of \eqref{eq:profile-single} determines the profile on each dyadic piece. Put $u=1+x$.
For $0\le x\le1/2$ it gives $\mathcal S(1+x)=\frac18+4\mathcal S(x)$. At $x=1/3$ the scaling relation turns this into
$\mathcal S(4/3)=\frac18+\frac14\mathcal S(4/3)$, hence
\begin{equation}\label{eq:profile-four-thirds}
 \mathcal S(4/3)=\tfrac16 .
\end{equation}
A function $h$ is H\"older continuous of exponent $0<\beta\le1$ on an interval if
$|h(u)-h(v)|\le K|u-v|^\beta$ for all $u,v$ in that interval and some constant $K$.

\begin{lemma}\label{lem:holder}
Let $\beta=\log_2(4/3)=2-\log_23$. The functions $\mathcal S$ and $\mathcal C$ are H\"older continuous of exponent $\beta$ on
$[0,2]$.
\end{lemma}

\begin{proof}
Let $(f_0,g_0)$ be the pair of affine functions of $\mathcal X$ and let $(f_n,g_n)=\mathcal T^n(f_0,g_0)$. On
$[1,2]$ the functions $f_0$ and $\mathcal S$ take their values in $[\frac18,\frac12]$, and $g_0$ and $\mathcal C$ in
$[\frac14,\frac12]$. The contraction therefore gives $|f_n-\mathcal S|\le\frac38(\frac34)^n$ and
$|g_n-\mathcal C|\le\frac38(\frac34)^n$ on $[1,2]$. Let $\Lambda_n$ be the larger of the Lipschitz constants of $f_n$
and $g_n$ on $[1,2]$. The extension of $f_n$ to $[0,1]$ has Lipschitz constant at most $\Lambda_n/2$, and the
extension of $g_n$ to $[0,2]$ has Lipschitz constant at most $\Lambda_n$. By \eqref{eq:profile-operator},
$\Lambda_{n+1}\le\frac32\Lambda_n$, and $\Lambda_0=\frac38$. Let $u$ and $v$ be in $[1,2]$ with $u\ne v$, and
choose $n\ge0$ with $2^{-n-1}<|u-v|\le2^{-n}$. Then
\[
 |\mathcal S(u)-\mathcal S(v)|\le\tfrac34\bigl(\tfrac34\bigr)^n+\tfrac38\bigl(\tfrac32\bigr)^n2^{-n}
 =\tfrac98\bigl(\tfrac34\bigr)^n\le\tfrac32\,|u-v|^{\beta},
\]
because $(\frac34)^n=2^{-n\beta}<(2|u-v|)^{\beta}$. The same bound holds for $\mathcal C$ on $[1,2]$. By scaling it holds
for $\mathcal S$ on each interval $[2^{-i},2^{1-i}]$ with $i\ge0$. Two points $x<y$ of $[0,2]$ with $y\le2x$ lie in the
same interval of this family or in two consecutive ones, which gives $\mathcal S(y)-\mathcal S(x)\le3(y-x)^{\beta}$. If $y>2x$
then $\mathcal S(y)-\mathcal S(x)\le \mathcal S(y)=O(y^2)=O((y-x)^2)$. Hence $\mathcal S$ is H\"older continuous of exponent $\beta$ on $[0,2]$, and
so is $\mathcal C$, which equals $2\mathcal S$ on $[0,1]$.
\end{proof}

\subsection{Convergence to the profile}

\begin{theorem}\label{thm:profile}
There is a constant $K$ such that, for every $k\ge1$,
\begin{equation}\label{eq:profile-convergence}
\begin{aligned}
 \max_{0\le n\le2^{k+1}}\left|\frac{\mathcal H(n)}{4^k}-\kappa \mathcal S\Bigl(\frac{n}{2^k}\Bigr)\right|&\le Kk^2\Bigl(\frac34\Bigr)^k,\\
 \max_{0\le t\le2^{k+1}}\left|\frac{\Gamma_{2^k}(t)}{4^k}-\kappa \mathcal C\Bigl(\frac{t}{2^k}\Bigr)\right|&\le Kk^2\Bigl(\frac34\Bigr)^k .
\end{aligned}
\end{equation}
Equivalently, the two cumulative counts differ from $\kappa4^k\mathcal S(n/2^k)$ and $\kappa4^k\mathcal C(t/2^k)$ by at most
$Kk^23^k$.
\end{theorem}

\begin{proof}
Let $e_k$ and $f_k$ be the two maxima in \eqref{eq:profile-convergence}. Let $\omega$ be a common modulus of
continuity of $\kappa \mathcal S$ and $\kappa \mathcal C$ on $[0,2]$. Fix $k$, put $N=2^k$, $q=N/2$ and
$\tau=\tau_k$, and assume $\tau+3\le q$. Put $\theta_k=|\mathcal M_q/N^2-\kappa/8|$. By
Theorem~\ref{thm:density} and the definition of $\mathcal M_q$ we have $\theta_k=O(k^3/2^k)$. We write $\eta_k$ for
quantities that are bounded by a fixed multiple of
$\theta_k+\omega\bigl((\tau+4)/N\bigr)+(\tau+2)^2/N^2+k^2/2^k$. By Lemma~\ref{lem:holder} we have $\omega(h)=O(h^{\beta})$, and $2^{-\beta k}=(\frac34)^k$. Since
$\tau=O(k^2)$, this gives $\eta_k=O\bigl(k^{2\beta}(\frac34)^k\bigr)$.

By Lemma~\ref{lem:degrees}, a diamond with center $m<q$ and size $p=\lb(m)$ is complete at round
$2m+2p+\tau\le N+\tau$. Hence $\mathcal H^{<q}(n)=\mathcal M_q$ for $n\ge N+\tau$. For $n\le N+2$ we have
$\mathcal H^{<q}(n)=\mathcal H(n)$, and $\mathcal H^{<q}\le \mathcal M_q$ everywhere. We also use the uniform crude bounds
$\mathcal H(n)=O((n+2)^2)$ and $\Gamma_r(t)=O((t+2)^2)$.
Every vertex born by time $t$ lies in a graph ball with $O((t+2)^2)$
vertices, independently of the cone size. The same bound holds for $a(n)$,
and \eqref{eq:an} gives $12\mathcal H(n-1)\le a(n)$ for $n\ge4$.
The remaining indices are checked directly.

For $n\le N$ the scaling of $\mathcal S$ shows that the difference in
\eqref{eq:profile-convergence} is one quarter of the corresponding difference at level $k-1$. It is at most
$e_{k-1}/4$. Let $n=N+J$ with $\tau\le J\le N$. By \eqref{eq:block-sum},
$\mathcal H(n)=\mathcal M_q+2\mathcal H(J)+\Gamma_q(J-2)-\epsilon$ with $\epsilon\in\{0,1\}$. By \eqref{eq:profile-operator} at the fixed point,
$\kappa \mathcal S(1+J/N)=\kappa/8+2\kappa \mathcal S(J/N)+\frac14\kappa \mathcal C(J/q)$. Subtracting and dividing by $N^2$, the term
with $\mathcal H(J)$ contributes at most $2e_{k-1}/4$, the term with $\Gamma_q$ contributes at most
$\frac14\bigl(f_{k-1}+\omega(2/q)\bigr)$, and the constants contribute at most $\theta_k+1/N^2$. The difference is therefore at most
$\frac12e_{k-1}+\frac14f_{k-1}+\eta_k$. Let $n=N+J$ with $0\le J<\tau$. Then
$\mathcal H(N)\le\mathcal H(n)\le \mathcal M_q+O((\tau+2)^2)$ by \eqref{eq:block-sum} and the crude bounds, while
$\kappa/8\le\kappa \mathcal S(n/N)\le\kappa/8+\omega(\tau/N)$. Since $\mathcal H(N)/N^2\ge\kappa/8-e_{k-1}/4$, the
difference is at most $\frac14e_{k-1}+\eta_k$. Altogether
\[
 e_k\le\tfrac12e_{k-1}+\tfrac14f_{k-1}+\eta_k .
\]

For the second error, put $M=q$, so that $N=2M$. Summing \eqref{eq:cone-coeff} gives, for $t\ge0$,
\[
 \Gamma_N(t)=E(t)+2\mathcal H^{<M}(t+1)+\Gamma_M(t-N)+\Gamma_{M-1}(t-N-1),
\]
with $|E(t)|\le3N$. For $t\le N-1$ the last two terms vanish and $\mathcal H^{<M}(t+1)=\mathcal H(t+1)$ with
$t+1\le N$. Since $\mathcal C(y)=2\mathcal S(y)$ on $[0,1]$, the difference in \eqref{eq:profile-convergence} is at most
$2e_{k-1}/4+2\omega(1/N)+3/N$. For $N+\tau\le t\le2N$ we have $\mathcal H^{<M}(t+1)=\mathcal M_M$, and
Lemma~\ref{lem:adjacent-uniform} replaces $\Gamma_{M-1}(t-N-1)$ by $\Gamma_M(t-N-1)$ with an error
$O(k^22^k)$. Since $\kappa \mathcal C(t/N)=\kappa/4+\frac12\kappa \mathcal C\bigl((t-N)/M\bigr)$ and $2\mathcal M_M/N^2$ differs from
$\kappa/4$ by $2\theta_k$, the difference is at most
$\frac14\bigl(2f_{k-1}+\omega(1/M)\bigr)+\eta_k$. For $N\le t<N+\tau$, monotonicity gives
$\Gamma_N(N-1)\le\Gamma_N(t)\le3N+2\mathcal M_M+O((\tau+2)^2)$, while $\kappa \mathcal C(t/N)$ lies between $\kappa/4$ and
$\kappa/4+\omega(\tau/N)$. The difference is at most $\frac12e_{k-1}+\eta_k$. Altogether
\[
 f_k\le\max\bigl(\tfrac12e_{k-1},\tfrac12f_{k-1}\bigr)+\eta_k .
\]

Let $v_k=\max(e_k,f_k)$. The two inequalities give $v_k\le\frac34v_{k-1}+\eta_k$ for all large $k$.
Multiplying by $(\frac43)^k$ and summing, we obtain $v_k=O\bigl(k^{1+2\beta}(\frac34)^k\bigr)$, and
$1+2\beta<2$.
\end{proof}

\subsection{The envelope and the optimal constants}

\begin{theorem}\label{thm:envelope}
For $1\le u\le2$,
\begin{equation}\label{eq:profile-envelope}
 \frac{3u^2}{32}\le \mathcal S(u)\le\frac{u^2}{8}.
\end{equation}
Equality holds on the left only at $u=4/3$, and on the right only at $u=1$ and $u=2$.
\end{theorem}

\begin{proof}
Let $\mathcal K$ be the set of pairs $(f,g)$ in $\mathcal X$ such that, on $[1,2]$,
\[
 \frac{3u^2}{32}\le f(u)\le\frac{u^2}{8},\qquad
 \frac14\le g(y)\le\frac12-\frac{(2-y)^2}{8}.
\]
It is closed and it contains the pair $f(u)=u^2/8$, $g(y)=y/4$. We show that $\mathcal T$ preserves
$\mathcal K$. Let $(f,g)\in\mathcal K$ and $(f^*,g^*)=\mathcal T(f,g)$. By scaling, the extension of $f$
satisfies the same two bounds on $[0,1]$, and the extension of $g$ satisfies $0\le g(y)\le y^2/4$ on $[0,1]$.
For $0\le x\le1/2$ we have $g(2x)=8f(x)$, and \eqref{eq:profile-operator} gives
\begin{equation}\label{eq:envelope-first}
\begin{aligned}
 f^*(1+x)-\frac{3(1+x)^2}{32}&=4\Bigl(f(x)-\frac{3x^2}{32}\Bigr)+\frac{(3x-1)^2}{32},\\
 \frac{(1+x)^2}{8}-f^*(1+x)&=4\Bigl(\frac{x^2}{8}-f(x)\Bigr)+\frac{x(2-3x)}{8}.
\end{aligned}
\end{equation}
For $1/2\le x\le1$ it gives
\begin{equation}\label{eq:envelope-second}
\begin{aligned}
 f^*(1+x)-\frac{3(1+x)^2}{32}&=2\Bigl(f(x)-\frac{3x^2}{32}\Bigr)+\frac14\Bigl(g(2x)-\frac14\Bigr)+\frac{3(1-x)^2}{32},\\
 \frac{(1+x)^2}{8}-f^*(1+x)&=2\Bigl(\frac{x^2}{8}-f(x)\Bigr)+\frac14\Bigl(\frac12-\frac{(1-x)^2}{2}-g(2x)\Bigr).
\end{aligned}
\end{equation}
All the terms on the right-hand sides are nonnegative, so $f^*$ satisfies the two bounds. For $0\le x\le1/2$
we have $0\le g(2x)\le x^2$, hence
\[
 \frac14\le g^*(1+x)\le\frac14+\frac{x^2}{2}\le\frac12-\frac{(1-x)^2}{8},
\]
the last inequality being $5x^2-2x-1\le0$. For $1/2\le x\le1$ the bounds on $g(2x)$ give
\[
 \frac38\le g^*(1+x)\le\frac12-\frac{(1-x)^2}{4}.
\]
Thus $\mathcal T$ preserves $\mathcal K$. The iterates of a point of $\mathcal K$ converge to the fixed point,
so $(\mathcal S,\mathcal C)\in\mathcal K$, which is \eqref{eq:profile-envelope}.

Take $f=f^*=\mathcal S$ and $g=g^*=\mathcal C$ in these identities. The lower bound gives $\mathcal S>0$ on $(0,\infty)$, hence
$\mathcal C>0$ on $(0,2]$ and $\mathcal C(y)>1/4$ for $1<y\le2$. In the first line of \eqref{eq:envelope-second} the middle
term is then positive for $x>1/2$, and the first line of \eqref{eq:envelope-first} vanishes only for $x=1/3$.
Together with \eqref{eq:profile-four-thirds} this proves the assertion on the left. For $0<x\le1/2$ the last
term of the second line of \eqref{eq:envelope-first} is positive. For $1/2\le x<1$ the two displayed bounds on
$g^*$ give $\mathcal C(2x)<\frac12-\frac12(1-x)^2$, so the second line of \eqref{eq:envelope-second} is positive. The
values $u=1$ and $u=2$ remain, where $\mathcal S(u)=u^2/8$.
\end{proof}

\begin{theorem}\label{thm:growth}
For $1\le u\le2$ let $\Psi(u)=12\kappa \mathcal S(u)/u^2$. Then, with $k=\lfloor\log_2n\rfloor$,
\begin{equation}\label{eq:growth}
 a(n)=n^2\,\Psi\Bigl(\frac{n}{2^k}\Bigr)+O\bigl(n^{\log_23}(\log n)^2\bigr).
\end{equation}
The function $\Psi$ is continuous, $\Psi(1)=\Psi(2)=3\kappa/2$, and
\begin{align*}
 \liminf_{n\to\infty}\frac{a(n)}{n^2}&=\frac{9\kappa}{8}=0.773905146264335750\ldots,\\
 \limsup_{n\to\infty}\frac{a(n)}{n^2}&=\frac{3\kappa}{2}=1.031873528352447667\ldots.
\end{align*}
The lower limit is approached along $n=\lfloor2^{k+2}/3\rfloor$ and the upper limit along $n=2^k$. In
particular $a(n)/n^2$ has no limit, and the quotient of the two constants is $3/4$.
\end{theorem}

\begin{proof}
Let $N=2^k\le n<2N$ and $u=n/N$. By Theorem~\ref{thm:profile} at the levels $k$ and $k+1$, each of
$\mathcal H(n-1)$, $\mathcal H(n)$, $\mathcal H(n+1)$, and $\mathcal H(n+2)$ differs from the corresponding value
of $\kappa N^2\mathcal S$ by $O(k^23^k)$. By Lemma~\ref{lem:holder} these four values of $\kappa N^2\mathcal S$ differ from
$\kappa N^2\mathcal S(u)$ by $O(N^{2-\beta})$, and $N^{2-\beta}=3^k$. Formula \eqref{eq:an} then gives
$a(n)=12\kappa N^2\mathcal S(u)+O(k^23^k)$, which is \eqref{eq:growth} because $3^k\le n^{\log_23}$. By Theorem~\ref{thm:envelope},
$9\kappa/8\le\Psi(u)\le3\kappa/2$ on $[1,2]$, with equality on the left at $u=4/3$ and on the right at $u=1$.
The two limits follow, the continuity of $\Psi$ being used for the sequence $n=\lfloor2^{k+2}/3\rfloor$.
\end{proof}

\begin{figure}[htbp]
 \centering
 \includegraphics[width=0.90\linewidth]{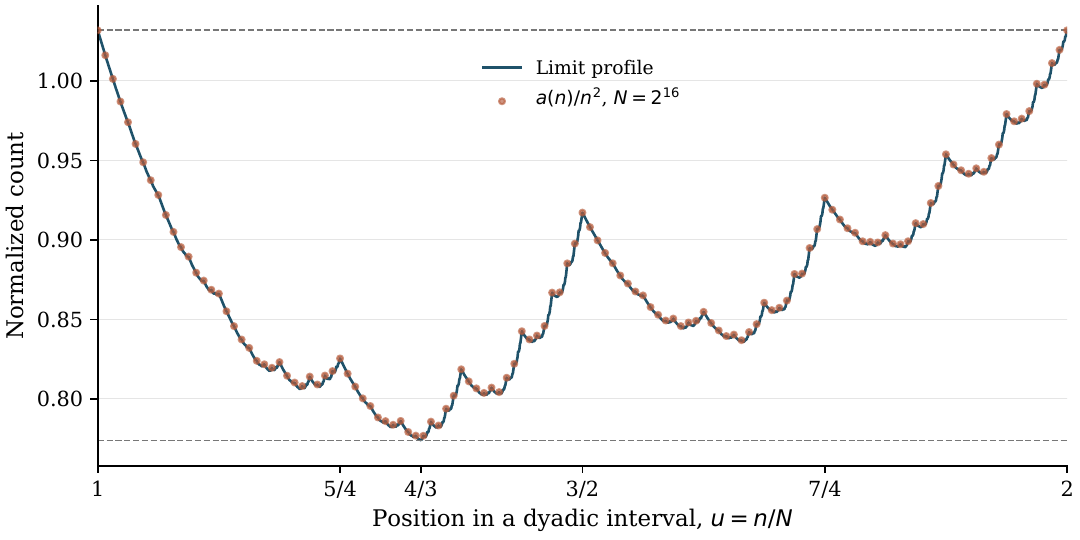}
 \caption{The function $\Psi(u)=12\kappa \mathcal S(u)/u^2$ on $[1,2]$, with values
 of $a(n)/n^2$ at $n=2^{16}u$. The dashed lines are the lower and upper limits of $a(n)/n^2$.}
 \label{fig:asymptotic}
\end{figure}

In logarithmic scale the main term is $n^2\Phi(\log_2n)$, where $\Phi(v)=\Psi(2^{v-\lfloor v\rfloor})$ is
continuous and has period~$1$. Figure~\ref{fig:asymptotic} compares the profile with the counts at scale $2^{16}$.

Theorem~\ref{thm:growth} is asymptotic. The next bound holds for every index.

\begin{proposition}
For every $n\ge0$ we have $a(n)\le(9n^2+3)/8$.
\end{proposition}

\begin{proof}
Each piece covers two edges of the lattice. By the proof of Lemma~\ref{lem:pieces}, a piece placed at a free
end covers the two sides that are not yet covered. Two free ends created at the same round have the same
color in the bipartite lattice and are not adjacent. Hence no edge is covered twice. The edges covered
during the first $n$ rounds join vertices at distance at most $n$ from the first apex. The sphere of radius
$d\ge1$ around a vertex of the honeycomb lattice has $3d$ vertices, and every edge joins two consecutive
spheres. The number $e_d$ of edges between the spheres of radius $d$ and $d+1$ therefore satisfies $e_0=3$ and
$e_d+e_{d-1}=9d$. The sum of $e_d$ over $d<n$ is $9n^2/4$ for even $n$ and $(9n^2+3)/4$ for odd $n$, and it
is at least $2a(n)$.
\end{proof}

\section{Binomial windows}\label{sec:windows}

The difference $X_k=C_{2^k}-C_{2^k-1}$ records the effect of adding the last layer of a dyadic cone.
The polynomials $\partial C_r$ and $E_c$ of Section~\ref{sec:cones} measure that effect from the first
possible birth in the added layer. Iterating \eqref{eq:transfer} decomposes the polynomial into
shifted groups of coefficients, which we call \emph{windows}, and a remainder computed on domains
of size $O(k)$.

\subsection{The finite decomposition}

The \emph{window polynomial} of depth $n\ge0$ is
\begin{equation}\label{eq:window-polynomial}
 \mathcal V_n(z)=\sum_{c=0}^{n}\binom nc E_c(z)-\bigl\lfloor2^{n-1}\bigr\rfloor .
\end{equation}
It has degree at most $J_n$. The first five are
\[
\begin{aligned}
 \mathcal V_0&=1,& \mathcal V_1&=1+z,& \mathcal V_2&=2+3z,\\
 \mathcal V_3&=4+6z+z^2+2z^3,&
 \mathcal V_4&=8+10z+5z^2+9z^3-z^5.
\end{aligned}
\]

\begin{theorem}\label{thm:windows}
Let $0\le D\le k$ and $q=2^{k-D}$, with $q\ge D$. Then
\begin{equation}\label{eq:window-decomposition}
\begin{aligned}
 \partial C_{2^k}(z)&=\sum_{n=0}^{D-1}z^{2^k-2^{k-n}}\,\mathcal V_n(z)+z^{2^k-q}R_{k,D}(z),\\
 R_{k,D}&=\sum_{c=0}^{D}\binom Dc\,\partial C_{q-c}-\bigl\lfloor2^{D-1}\bigr\rfloor .
\end{aligned}
\end{equation}
For $k\ge1$, choose $q=2^{\lceil\log_2(k+1)\rceil}$ and $D=k-\log_2q$. Then every kernel and every remainder term is computed on cones and triangles of size at
most $2k+1$. Moreover
\begin{equation}\label{eq:polynomial-support}
 \#\operatorname{supp}(X_k)=O(k^3).
\end{equation}
\end{theorem}

\begin{proof}
We prove by induction on $N$, for $0\le N\le D$, that
\[
 \partial C_{2^k}=\sum_{n=0}^{N-1}z^{2^k-2^{k-n}}\Bigl(\sum_{c=0}^{n}\binom nc E_c-2^n z^{2^{k-n-1}}\Bigr)
 +z^{2^k-2^{k-N}}\sum_{c=0}^{N}\binom Nc\,\partial C_{2^{k-N}-c}.
\]
For $N=0$ the identity is immediate. Let $N<D$ and put $M=2^{k-N-1}$. Every index $c$ of the last sum
satisfies $c\le N\le D-1<q\le M$. Proposition~\ref{prop:transfer} therefore applies to
$\partial C_{2M-c}$. The terms $\partial C_{M-c}$ and $\partial C_{M-c-1}$ are collected with Pascal's rule
and carry the factor $z^{2^k-2M+M}=z^{2^k-2^{k-N-1}}$. The constants add up to $-2^N$ with the same factor.
This is the statement for $N+1$. In the statement for $N=D$, the monomial $-2^nz^{2^k-2^{k-n-1}}$ is the
constant term $-2^n$ of the window of depth $n+1$ when $n+1\le D-1$, and it joins the last sum when
$n=D-1$. This gives \eqref{eq:window-decomposition}.

For the stated choice we have $k+1\le q<2(k+1)$ and $D=k-\lceil\log_2(k+1)\rceil<q$. The remainder uses the
cones of size at most $q$. The windows use the kernels $E_c$ with $c<D$, which are computed in triangles of
size at most $2k$. By Lemma~\ref{lem:kernel-support} the remainder has $O(q^2)$ possible exponents and the
window of depth $n$ has at most $J_n+1$. Summing over $n<D$ gives \eqref{eq:polynomial-support}.
\end{proof}

\begin{figure}[htbp]
 \centering
 \includegraphics[width=0.95\linewidth]{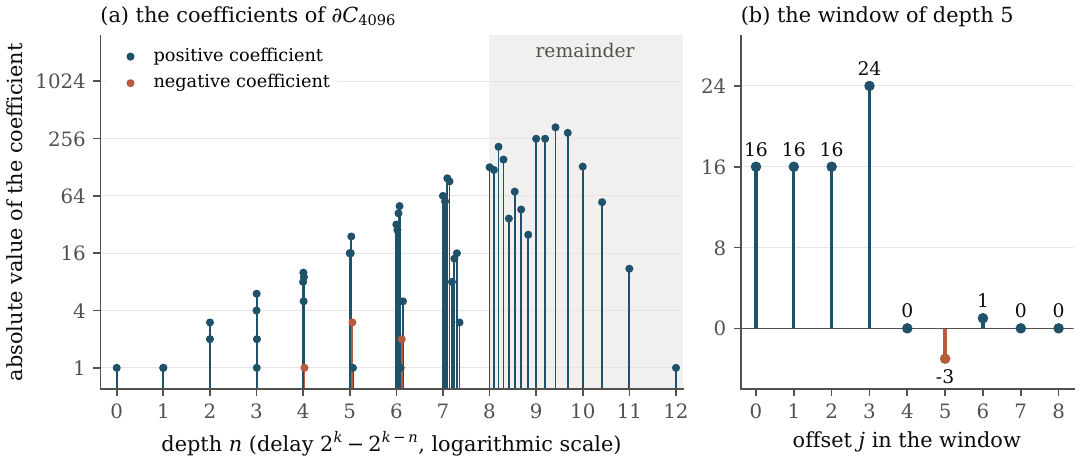}
 \caption{(a) The $51$ nonzero coefficients of $\partial C_{4096}$. The abscissa is $-\log_2(1-t/4096)$ for
 the delay $t$, so that the window of depth $n$ starts at $n$. The shaded part is the remainder of
 \eqref{eq:window-decomposition} for $q=16$. (b) The window of depth $5$, which is $\mathcal V_5$.}
 \label{fig:windows}
\end{figure}

The window of depth $n$ starts at the delay $2^k-2^{k-n}$. The windows therefore accumulate toward the delay
$2^k$, and the distance between two consecutive windows is halved at each step. Figure~\ref{fig:windows}
shows the case $k=12$.

\begin{corollary}
Let $k=a+n$ with $a\ge1$ and $n\ge0$. If $2^{a-1}>\max(n,J_n)$, then for $0\le j<2^{a-1}$
\[
 [z^{\,2^k-2^a+j}]\,\partial C_{2^k}(z)=[z^j]\,\mathcal V_n(z).
\]
In particular this coefficient depends only on $n$ and $j$.
\end{corollary}

\begin{proof}
Use Theorem~\ref{thm:windows} with $D=n+1$, for which $q=2^{a-1}$. The remainder starts at the delay
$2^k-2^{a-1}$. A window of depth $n'<n$ ends before the delay $2^k-2^{k-n'}+2^{a-1}$, which is at most
$2^k-2^a$. The window of depth $n$ starts at $2^k-2^a$. Hence it is the only term that contributes to the
delays from $2^k-2^a$ to $2^k-2^{a-1}-1$.
\end{proof}

\subsection{The first rounds after a cut}

In a growth, the \emph{first-touch data} of a vertex other than the seed are the time at which it first
receives a proposal and the number of proposals that it receives at that time. We assign the seed the initial
data $(0,1)$. A non-seed vertex that never receives a proposal has no data. With this convention, a vertex is
born at time $t$ exactly when its data are $(t,1)$. Along an edge of the lattice either $g$ changes by one
and $y$ is unchanged, or $y$ changes by one and $g$ is unchanged.

\begin{lemma}\label{lem:causal}
Let $0\le c<M$ and $t_0=2M-c-1$. Consider the growths on $\Omega_{M,c}$ and on $\Omega_{M,c+1}$. Let $v$ be a
vertex of $\Omega_{M,c+1}$ and write $g(v)=t_0-\iota$, with $\iota\ge1$. If the first-touch data of $v$ in the
two growths differ, then $v$ receives no proposal before time $t_0+\iota+y$ in either growth, where $y$ is the
height of $v$.

Consequently, let $j\ge0$ and let $v$ be a vertex of $\Omega_{M,c}$ that is born at time $t_0+j$ in exactly
one of the two growths. If $g(v)<t_0$ then $\iota+y\le j$. If $g(v)=t_0$ then $y\le j$.
\end{lemma}

\begin{proof}
Let $T$ be the least time at which $v$ receives a proposal in one of the two growths. We prove
$T\ge t_0+\iota+y$ by induction on $T$. The seed has the same data in both growths, so $T\ge1$. Before time
$T$ the vertex $v$ is intact in both growths, and the numbers of proposals that it receives at time $T$
differ. These proposals are sent by the neighbors of $v$ born at time $T-1$. Hence some neighbor $w$ of $v$
is born at time $T-1$ in exactly one of the two growths. Suppose first that $w\notin\Omega_{M,c+1}$. Then
$g(w)=t_0$ and $g(v)\le t_0-1$, so $g$ changes along the edge. Therefore $\iota=1$ and $w$ has height $y$. The vertex $w$ is born at time $T-1$ in the growth on $\Omega_{M,c}$, which is at least its
earliest time $t_0+y$. Thus $T\ge t_0+1+y$. Suppose now that $w\in\Omega_{M,c+1}$, and write
$g(w)=t_0-\iota'$ and $y'$ for its height. The data of $w$ in the two growths differ, and $w$
receives a proposal at time $T-1$ in one of them. By the induction hypothesis $T-1\ge t_0+\iota'+y'$. Along
the edge the sum $\iota+y$ changes by one, so $\iota'+y'\ge\iota+y-1$ and $T\ge t_0+\iota+y$.

For the second assertion, suppose that $g(v)<t_0$. The data of $v$ differ and $v$ receives a proposal at time
$t_0+j$ in one of the growths, so $t_0+j\ge t_0+\iota+y$. Suppose that $g(v)=t_0$. Then $v$ is born at time
$t_0+j$ in the growth on $\Omega_{M,c}$, and this time is at least its earliest time $t_0+y$.
\end{proof}

By \eqref{eq:kernel}, the coefficient $e_j(c)$ is the number of vertices born at time $t_0+j$ in the growth on
$\Omega_{M,c}$ minus the same number for $\Omega_{M,c+1}$. The lemma confines the vertices that contribute to
this difference to a triangle of side $j$ next to the last layer. For $j\le1$, this determines the coefficients.

\begin{proposition}\label{prop:first-window-columns}
For every $c\ge0$,
\[
 e_0(c)=1,\qquad e_1(c)=\mathbf1_{c\equiv1,2\pmod4}.
\]
Consequently
\[
 [z^0]\mathcal V_n=
 \begin{cases}1,&n=0,\\2^{n-1},&n\ge1,\end{cases}
 \qquad
 [z]\mathcal V_n=
 \sum_{\substack{0\le c\le n\\c\equiv1,2\ (4)}}\binom nc
 =\frac{2^n-\operatorname{Re}(1+i)^{n+1}}2.
\]
\end{proposition}

\begin{proof}
The first formula is the constant term of $E_c$ in Lemma~\ref{lem:kernel-support}. Let $M$ be dyadic with
$M>c$ and let $v$ be a vertex born at time $t_0+1$ in exactly one of the two growths of
Lemma~\ref{lem:causal}. If $g(v)<t_0$ then $\iota=1$ and $y=0$. By Lemma~\ref{lem:pascal} such a vertex of
the axis is born at its earliest time $t_0-1$ in both growths, which is a contradiction. Hence $g(v)=t_0$
and $y\le1$. The vertex of the axis with $g=t_0$ is born at time $t_0$. Therefore $v$ is the vertex with
$g=t_0$ and $y=1$, and $e_1(c)$ equals $1$ if this vertex belongs to $\Omega_{M,c}$ and is born at its
earliest time $t_0+1$, and $0$ otherwise. This vertex has $2x+p=2M-c-2$. By Lemma~\ref{lem:pascal} it is born
at its earliest time exactly when $x$ is even. For $c=0$ it is $A(M-1,1)$, which is outside the triangle.
Let $c\ge1$, so that $M$ is even. If $c$ is even the vertex is $A(M-1-c/2,1)$ and $x$ is even exactly when
$c\equiv2\pmod4$. If $c$ is odd the vertex is $B(M-(c+3)/2,1)$ and $x$ is even exactly when
$c\equiv1\pmod4$. In both cases $x+1<M$, so the vertex belongs to $\Omega_{M,c}$. This proves the second
formula. Identity \eqref{eq:window-polynomial} gives the binomial sums, and the fourth roots of unity give the
displayed expression for the second one.
\end{proof}

\subsection{Periodic kernels}

The next lemma compares two truncated triangles from two columns whose abscissas are multiples of a power of
two.

\begin{lemma}\label{lem:young}
Let $P$ be a power of two. Let $\Omega=\Omega_{M,c}$ and $\Omega'=\Omega_{M',c'}$, and let $k$ and $k'$ be
multiples of $P$ with $0\le k<M$ and $0\le k'<M'$. Assume that $2M-c-2k=2M'-c'-2k'$. Put
$m=\min(M-k,M'-k')$, and let $R$ be the set of triples $(p,x,y)$ with $x\ge0$, $y\ge0$, $x+y<m$ and
$(p,k+x,y)\in\Omega$. A triple $(p,x,y)$ with $x\ge0$, $y\ge0$ and $x+y<m$ belongs to $R$ if and only if
$(p,k'+x,y)\in\Omega'$. Moreover, for every $(p,x,y)\in R$, every integer $s$ with
\[
 s<2x+p+2P\qquad\text{and}\qquad s\le2m,
\]
and every $\mu\ge1$, the vertex $(p,k+x,y)$ has first-touch data $(2k+s,\mu)$ in the growth on $\Omega$ if and
only if the vertex $(p,k'+x,y)$ has first-touch data $(2k'+s,\mu)$ in the growth on $\Omega'$.
\end{lemma}

\begin{proof}
For $x+y<m$ the vertex $(p,k+x,y)$ satisfies the first constraint of $\Omega$, and it satisfies the second
one exactly when $2x+y+p<2M-c-2k$. The same holds for $\Omega'$, which proves the first assertion. We call
$s$ the relative time, and we say that a triple $\mathbf v=(p,x,y)$ of $R$ has the event $(s,\mu)$ in
$\Omega$ when $(p,k+x,y)$ has first-touch data $(2k+s,\mu)$, and similarly in $\Omega'$. Put
$h(\mathbf v)=2x+p+2P$. Along an edge $h$ changes by one. The earliest time of $(p,k+x,y)$ is
$2k+2(x+y)+p$, so no triple of $R$ has an event at a negative relative time. We prove the statement by
induction on $s$. Let $\mathbf v\in R$ and let $s$ satisfy $s<h(\mathbf v)$ and $s\le2m$, and assume the
statement for the smaller relative times.

Both truncated triangles are closed under predecessors for the earliest-time potential $2(x+y)+p$, so
Lemma~\ref{lem:pascal} applies to them. Suppose first that $\mathbf v=(0,0,y)$. If $y<P$ then
$k\bitand y=0$ and $k'\bitand y=0$, since $k$ and $k'$ are multiples of $P$. By Lemma~\ref{lem:pascal} the
vertices $A(k,y)$ and $A(k',y)$ are born at their earliest times, so $\mathbf v$ has the event $(2y,1)$ in
both domains and no other event. If $y\ge P$ the earliest relative time of $\mathbf v$ is
$2y\ge2P=h(\mathbf v)>s$, and $\mathbf v$ has no event at the relative time $s$ in either domain.

Suppose now that $\mathbf v$ is not of this form. Then every neighbor of $(p,k+x,y)$ has an $x$-coordinate
at least $k$, and the same holds with $k'$. By the induction hypothesis, applied to $\mathbf v$ and to the
relative times smaller than $s$, the triple $\mathbf v$ is intact before the relative time $s$ in $\Omega$
exactly when it is so in $\Omega'$. The proposals that reach it at the relative time $s$ are sent by its
neighbors born at the relative time $s-1$. A neighbor $\mathbf w\in R$ satisfies
$h(\mathbf w)\ge h(\mathbf v)-1>s-1$ and $s-1\le2m$, so by the induction hypothesis it is born at the
relative time $s-1$ in $\Omega$ exactly when it is so in $\Omega'$. A neighbor $(p',x',y')$ with $x'\ge0$,
$y'\ge0$ and $x'+y'\ge m$ has an earliest relative time at least $2m>s-1$, so it is not born at the relative
time $s-1$ in either domain. The remaining neighbors have a negative coordinate or lie outside both domains,
by the first assertion. Hence $\mathbf v$ receives the same number of proposals at the relative time $s$ in
both domains, which completes the induction.
\end{proof}

\begin{theorem}\label{thm:kernel-columns}
Let $j\ge1$ and let $P$ be a power of two with $P>j$. Then
\[
 e_j(c+2P)=e_j(c)\qquad\text{for every } c\ge j-1 .
\]
\end{theorem}

\begin{proof}
Let $c\ge j-1$ and $\hat c=c+2P$. Let $M$ be a dyadic integer with $M>\hat c$. Then $P$ divides $M$. Put
\[
 t_0=2M-c-1,\qquad \hat t_0=t_0-2P,\qquad m=M-P,\qquad s=\hat t_0+j=2m-(c+1-j).
\]
Thus $s\le2m$, with equality exactly when $c=j-1$. For a vertex $v$ of $\Omega_{M,c}$ put
\[
 \sigma(v)=\mathbf 1\bigl[v \text{ is born at time } t_0+j \text{ in } \Omega_{M,c}\bigr]
          -\mathbf 1\bigl[v \text{ is born at time } t_0+j \text{ in } \Omega_{M,c+1}\bigr],
\]
where the second term is zero when $v\notin\Omega_{M,c+1}$. Define $\hat\sigma(v)$ in the same way for a
vertex $v$ of $\Omega_{M,\hat c}$, with the cuts $\hat c$ and $\hat c+1$ and the time $\hat t_0+j=s$. By
\eqref{eq:kernel},
\begin{equation}\label{eq:two-sums}
 e_j(c)=\sum_{v\in\Omega_{M,c}}\sigma(v),\qquad e_j(\hat c)=\sum_{v\in\Omega_{M,\hat c}}\hat\sigma(v).
\end{equation}
We say that a triple $(p,x,y)$ with $x\ge0$ and $y\ge0$ is \emph{sensitive} when, with $\gamma=2x+y+p$, either
$\gamma=\hat t_0-\iota$ with $\iota\ge1$ and $\iota+y\le j$, or $\gamma=\hat t_0$ and $y\le j$. A sensitive
triple satisfies $y\le j$ and
\begin{equation}\label{eq:sensitive-bounds}
 2x+p\ge s-2j,\qquad 2(x+y)+p\le s .
\end{equation}
In the second inequality, equality requires $\gamma=\hat t_0$ and $y=j$.

By Lemma~\ref{lem:causal}, applied with the cut $\hat c$, every vertex with $\hat\sigma\ne0$ is a sensitive
triple. By \eqref{eq:sensitive-bounds} it satisfies $x+y\le m$, and $x+y=m$ is possible only if $c=j-1$,
$p=0$ and $y=j$. Let $v=(p,x,y)$ be a vertex with $\sigma(v)\ne0$. By Lemma~\ref{lem:causal}, applied with the
cut $c$, we have $y\le j$ and $g(v)\ge t_0-j$, hence $2x\ge t_0-2j-1$. Since $t_0\ge2\hat c+1-c\ge4P+j$, this
gives $x\ge P$. Moreover $(p,x-P,y)$ is a sensitive triple, because $g(v)-t_0=(g(v)-2P)-\hat t_0$, and
$(x-P)+y<m$ because $v\in\Omega_{M,0}$.

Let $\mathcal R$ be the set of sensitive triples with $x+y<m$, and let $(p,x,y)\in\mathcal R$. The vertex
$(p,P+x,y)$ belongs to $\Omega_{M,c}$ and the vertex $(p,x,y)$ belongs to $\Omega_{M,\hat c}$. The first
belongs to $\Omega_{M,c+1}$ exactly when the second belongs to $\Omega_{M,\hat c+1}$, that is when
$\gamma<\hat t_0$. By \eqref{eq:sensitive-bounds}, $2x+p+2P\ge s-2j+2P>s$, and $s\le2m$. We apply
Lemma~\ref{lem:young} with $k=P$ and $k'=0$, for which $\min(M-k,M-k')=m$. With the domains $\Omega_{M,c}$
and $\Omega_{M,\hat c}$, whose relative cuts are both $\hat t_0+1$, it shows that $(p,P+x,y)$ is born at time
$2P+s=t_0+j$ in the growth on $\Omega_{M,c}$ if and only if $(p,x,y)$ is born at time $s$ in the growth on
$\Omega_{M,\hat c}$. With the domains $\Omega_{M,c+1}$ and $\Omega_{M,\hat c+1}$, whose relative cuts are both
$\hat t_0$, it shows the same equivalence for these two growths when $\gamma<\hat t_0$. Therefore
$\sigma(p,P+x,y)=\hat\sigma(p,x,y)$.

By the two preceding paragraphs, the first sum of \eqref{eq:two-sums} is the sum of $\hat\sigma$ over
$\mathcal R$, and the second sum is the sum of $\hat\sigma$ over $\mathcal R$ and over the sensitive triples
with $x+y=m$. If $c>j-1$ there is no such triple, so $e_j(\hat c)=e_j(c)$. Let $c=j-1$. Then $s=2m$ and
$\hat t_0=2m-j$, and the only sensitive triple with $x+y=m$ is the vertex $A(m-j,j)$. It belongs to
$\Omega_{M,\hat c}$ and not to $\Omega_{M,\hat c+1}$, and $2m$ is its earliest time. By
Lemma~\ref{lem:pascal} it is born at time $2m$ in the growth on $\Omega_{M,\hat c}$ exactly when
$(m-j)\bitand j=0$. This fails. Indeed $1\le j<P$ and $P$ divides $m$, so if $2^\ell$ is the largest power of
two that divides $j$, then $j$ and $m-j$ are both odd multiples of $2^\ell$. Hence $\hat\sigma$ vanishes at
this vertex and $e_j(\hat c)=e_j(c)$.
\end{proof}

The proof uses no bound on the time of the last birth in a triangle. For $c=j-1$ the coefficient $e_j(c)$ is
read at time $2M$, and Lemma~\ref{lem:young} is used at its last admissible relative time $s=2m$.
The rows $c\le j-2$ are not covered. They concern births after time $2M$.

\subsection{The columns of the windows}

For $j\ge1$ let $p_j=2^{\lfloor\log_2j\rfloor+2}$. It is twice the least power of two larger than $j$, so
Theorem~\ref{thm:kernel-columns} gives $e_j(c+p_j)=e_j(c)$ for $c\ge j-1$. Define
\begin{equation}\label{eq:pi-eta}
\begin{aligned}
 \pi_j(r)&=e_j\bigl(j-1+((r-j+1)\bmod p_j)\bigr) &&(0\le r<p_j),\\
 \eta_j(c)&=e_j(c)-\pi_j(c\bmod p_j) &&(0\le c\le j-2).
\end{aligned}
\end{equation}
The list $\pi_j$ is one period of the column $j$, and $\eta_j$ measures the defect of periodicity on the first
rows. The list $\eta_j$ is empty for $j=1$. For $n\ge0$ let $v_j(n)=[z^j]\mathcal V_n(z)$.

\begin{proposition}\label{prop:window-columns}
Let $j\ge1$. For every $n\ge0$,
\begin{equation}\label{eq:window-column}
 v_j(n)=\sum_{r=0}^{p_j-1}\pi_j(r)\sum_{\substack{0\le c\le n\\ c\equiv r\ (\mathrm{mod}\ p_j)}}\binom nc
 +\vartheta_j(n),\qquad
 \vartheta_j(n)=\sum_{c=0}^{j-2}\eta_j(c)\binom nc .
\end{equation}
The function $\vartheta_j$ is a polynomial in $n$ of degree at most $j-2$, and $\vartheta_1=0$.
\end{proposition}

\begin{proof}
By \eqref{eq:window-polynomial}, $v_j(n)=\sum_{c=0}^{n}\binom nc e_j(c)$. We have
$e_j(c)=\pi_j(c\bmod p_j)$ for $c\ge j-1$ by Theorem~\ref{thm:kernel-columns}, and
$e_j(c)=\pi_j(c\bmod p_j)+\eta_j(c)$ for $c\le j-2$ by \eqref{eq:pi-eta}.
\end{proof}

\begin{corollary}\label{cor:depth-recurrence}
Let $\nabla$ be the difference operator, $(\nabla f)(n)=f(n+1)-f(n)$. For every $j\ge1$,
\[
 \nabla^{\,j-1}\bigl(\nabla^{\,p_j}-1\bigr)v_j=0 .
\]
\end{corollary}

\begin{proof}
For a residue $r$ modulo $p_j$ let $f_r(n)$ be the inner sum of \eqref{eq:window-column}. Pascal's rule gives
$\nabla f_r=f_{r-1}$, hence $\nabla^{p_j}f_r=f_r$. It also gives $\nabla\binom nc=\binom n{c-1}$, so
$\nabla^{j-1}$ annihilates $\vartheta_j$. The two operators commute.
\end{proof}

These are recurrences for each fixed column, and their order grows with $j$. The number of sequences needed as $j$ varies is not bounded here.

For $1\le j\le63$, a finite calculation gives $\eta_j=0$ and proves that $p_j$ is the least period.
Appendix~\ref{app:columns} gives the procedure, the first seven lists, and the certificate.
The general theorem requires only eventual periodicity and allows a nonzero polynomial correction.

\section{Open questions}

The optimal error term in \eqref{eq:growth} and the best H\"older exponent
of the profile are unknown. It is also open whether
$9\kappa/8<a(n)/n^2<3\kappa/2$ holds for every $n\ge3$.

The identity $\deg P_q=2q-3$ remains open beyond the range certified in
Appendix~\ref{app:bounds}. By Lemma~\ref{lem:degrees}, it would follow from
$\deg C_r=2r-2$ for every direct cone. A bound $J_n=O(n)$ for the kernel degrees would improve
\eqref{eq:polynomial-support} to $\#\operatorname{supp}(X_k)=O(k^2)$.

Theorem~\ref{thm:kernel-columns} gives eventual periodicity for every kernel column.
It is unknown whether all prefix corrections $\eta_j$ vanish, or equivalently whether, for every
dyadic $P$, the coefficients of $E_c$ of degree less than $P$ are determined by $c$ modulo $2P$.
A rule for the period lists $\pi_j$ is also unknown.

The series of Theorem~\ref{thm:ternary} computes $\kappa$ with a geometric
error bound. Its rationality, algebraicity, or transcendence is unknown.
A Mellin analysis of the periodic function $\Phi$, along the lines of
Flajolet, Gourdon, and Dumas~\cite{FGD}, may describe its Fourier coefficients.
The analytic continuation of the generating function is another question.
The denominators $1-z^{4q}$ alone do not exclude cancellation in the
infinite sum.

It is not known whether finitely many auxiliary numerical sequences give $a(n)$ by recurrences
of a fixed form, or whether the direct cones can be eliminated. Another question is whether all
subsequences $(a(2^e n+r))_{n\ge0}$, with $e\ge0$ and $0\le r<2^e$, are integer linear combinations
of finitely many fixed integer sequences, as in Allouche and Shallit~\cite{AS}.
No equivalence between these questions is asserted. It is also unknown whether a nontrivial linear
equation with polynomial coefficients relates $A(z)$ to finitely many of the series
$A(z^2),A(z^4),\ldots$, or whether the complexity bound can be improved.

The geometric method may apply to other toothpick constructions.
It requires separating boundaries and a bound on the unreduced regions,
in addition to the earliest-time Pascal pattern.

\appendix
\section{Explicit bounds and finite degree checks}\label{app:bounds}

The estimates used in the text can be made effective without an assumption about the last birth time.
The sizes of the cone and triangle give
\[
 0\le C_r(1)\le\frac{r(r+1)}2,\qquad 0\le L_{m,c}(1)\le m(m+1).
\]
For $c\ge1$, the least dyadic integer $m>c$ satisfies $m\le2c$. Consequently
\[
 \|E_c\|\le2m(m+1)\le12c^2,
\]
while $E_0=1$. These inequalities justify the following constants in the proof of
Lemma~\ref{lem:adjacent-uniform}. An internal node has forcing norm at most $29(k+1)^2$.
At a leaf, the two cones have index smaller than $4(k+1)$, so their difference has norm less than
$16(k+1)^2$. Its value at $1$ has absolute value at most $10(k+1)^2$, since both masses are nonnegative.
Counting the leaves and internal nodes gives
\[
 \|C_{2^k}-C_{2^k-1}\|\le45(k+1)^2 2^k,\qquad
 |d_k|\le40(k+1)^2 2^k.
\]
For a tree starting from $Y(a,j)$, the same argument gives
$|Y(a,j)|\le40(a+j+1)^2 2^a$.
Since $d_1=2$, identity \eqref{eq:mass-increments} yields
$|\delta_k|\le81(k+1)^2 2^k$. Therefore
\begin{equation}
 \left|\kappa-\frac{\mu_{2^k}}{4^k}\right|
 \le\frac{81}{4}\sum_{s\ge k}\frac{(s+1)^2}{2^s}
 =\frac{81(k^2+4k+6)}{2^{k+1}}
 \le\frac{42(k^2+4k+6)}{2^k}.
\end{equation}
The value $\mu_{2^{32}}=12\,689\,804\,595\,979\,239\,659$ and this bound imply $\kappa>0$.
At $k=128$, the bound is less than $2.087\cdot10^{-33}$.
This mass, computed by both programs, and rational interval arithmetic certify the thirty
decimal digits in Theorem~\ref{thm:density}. The integer and the interval endpoints are included with the code.

For the degree excess of Lemma~\ref{lem:degrees}, the domain bound gives $\tau_k\le5k^2$.
A direct cone of index $r\ge4$ has $r=2^a-c$, $c>2^{a-2}$, and rank $a+c$.
If this rank is at most $266$, then $a\le9$ and $r\le383$.
The computation of the cones through index $383$ gives $\deg C_r=2r-2$, so
$\tau_k=0$ for $k\le266$. Lemma~\ref{lem:degrees} then gives $\deg P_q=2q-3$ for dyadic $q\le2^{266}$.
Let $N=2q$ with $q\ge2$ dyadic. The term $F^{<q}(N+j)$ in \eqref{eq:block-rule} counts births of diamonds with center below $q$ at an age
larger than their size allows when $\deg P_p=2p-3$. When this degree condition holds for every $p<q$, the
term vanishes and Lemmas~\ref{lem:block} and~\ref{lem:cone-coeff} read
\begin{equation}\label{eq:exact-rule}
\begin{aligned}
 F(N+j)&=2F(j)+c_{N/2}(j-2) &&(3\le j\le N),\\
 c_N(t)&=2+2F(t+1) &&(2\le t\le N-1),\\
 c_N(N+u)&=c_{N/2}(u)+c_{N/2-1}(u-1)-\mathbf 1_{1\le u\le N/2-1} &&(u\ge0).
\end{aligned}
\end{equation}
This gives the following formula at dyadic times.

\begin{proposition}
Let $N=2^k$ with $3\le k\le268$. Then $F(N)=N/4-1$ and
\begin{equation}\label{eq:dyadic-count}
 a(N)=12\,\mathcal M_{N/2}+5N-7 .
\end{equation}
\end{proposition}

\begin{proof}
By the finite verification above and Lemma~\ref{lem:degrees}, $\deg P_p=2p-3$ and $\deg C_p=2p-2$ for every dyadic $p\le N/4$, and
$\deg C_{p-1}=2p-4$ for every dyadic $p$ with $2\le p\le N/8$. In particular the identities \eqref{eq:exact-rule}
hold at every dyadic scale from $4$ to $N/2$. Let $r$ be dyadic with $4\le r\le N/4$. The third line of \eqref{eq:exact-rule} at scale $r$ with $u=r-2$ gives
$c_r(2r-2)=c_{r/2}(r-2)$, because the coefficient of $z^{r-3}$ in $C_{r/2-1}$ is zero. Hence $c_r(2r-2)=1$ by induction from
$c_1(0)=c_2(2)=1$. The first line at scale $N/2$ with $j=N/2$ then gives $F(N)=2F(N/2)+1$. Since $F(4)=0$ we obtain
$F(N)=N/4-1$. By the degree condition every diamond with center below $N/2$ is complete at round $N$, and no
other diamond appears before round $N+3$. Thus $\mathcal H(n)=\mathcal M_{N/2}$ for $n\in\{N,N+1,N+2\}$ and
$\mathcal H(N-1)=\mathcal M_{N/2}-F(N)$. Formula \eqref{eq:an} gives \eqref{eq:dyadic-count}.
\end{proof}

\section{Certification of the first kernel columns}\label{app:columns}

Theorem~\ref{thm:kernel-columns} holds for every $j$. The next proposition treats the rows $c\le j-2$ by a finite computation for $j\le63$.

\begin{proposition}
For $1\le j\le63$ the list $\eta_j$ vanishes. Hence $e_j(c)=\pi_j(c\bmod p_j)$ for every $c\ge0$,
$\vartheta_j=0$, and $\bigl(\nabla^{p_j}-1\bigr)v_j=0$. The lists $\pi_j$ for $j\le7$ are
\[
\begin{array}{lll}
 j=1 & p_1=4 & 0,\,1,\,1,\,0\\
 j=2 & p_2=8 & 0,\,0,\,0,\,1,\,1,\,1,\,1,\,0\\
 j=3 & p_3=8 & 0,\,0,\,0,\,2,\,1,\,-1,\,1,\,0\\
 j=4 & p_4=16 & 0,\,0,\,0,\,0,\,0,\,0,\,1,\,1,\,1,\,1,\,1,\,1,\,1,\,-1,\,2,\,0\\
 j=5 & p_5=16 & 0,\,0,\,0,\,0,\,-1,\,2,\,1,\,0,\,1,\,2,\,1,\,-1,\,0,\,1,\,2,\,0\\
 j=6 & p_6=16 & 0,\,0,\,0,\,0,\,0,\,1,\,-1,\,2,\,2,\,1,\,1,\,-1,\,-1,\,2,\,1,\,0\\
 j=7 & p_7=16 & 0,\,0,\,0,\,0,\,0,\,0,\,0,\,3,\,1,\,-1,\,2,\,-1,\,0,\,1,\,1,\,0 .
\end{array}
\]
\end{proposition}

\begin{proof}
For $1\le j\le J$, formula \eqref{eq:pi-eta} uses $e_j(c)$ only for
$0\le c\le j-2+p_j\le J-2+p_J$. For each such $c$, choose the least dyadic $m>c$ and simulate the
two finite domains $\Omega_{m,c}$ and $\Omega_{m,c+1}$. Form $E_c$ from their complete birth polynomials,
read the required coefficients, and apply \eqref{eq:pi-eta}. Each simulation stops by
Remark~\ref{rem:finite-domain}. For $J=63$, the computation uses $c\le189$ and $m\le256$.
It gives $\eta_j=0$ for every $j\le63$. The other assertions follow from Theorem~\ref{thm:kernel-columns},
Proposition~\ref{prop:window-columns}, and Corollary~\ref{cor:depth-recurrence}.
\end{proof}

The list $\pi_1$ agrees with Proposition~\ref{prop:first-window-columns}. The computation was repeated with a
second program, written separately. In the computed range $c<260$ every kernel satisfies $\deg E_c\le c+1$,
which means that the triangles $\Omega_{m,c}$ and $\Omega_{m,c+1}$ have the same numbers of births at every
time after $2m$. This is a finite observation and not a general statement.

In this section a period is any integer $T$ for which the stated periodicity holds. For $1\le j\le63$ the
computed lists also show that $\pi_j$ does not have the period $p_j/2$. Since the least period of a sequence
divides each of its periods, $p_j$ is the least period of $e_j$ for these columns. We make no claim of
minimality for $j\ge64$.

\clearpage
\section*{Use of AI}

The author used OpenAI's ChatGPT 6 Astra to explore the V-toothpick sequence.
The first prompt, posed partly in jest, was: ``Find a formula or an asymptotic formula for $a(n)$.''
The AI supplied neither a usable answer nor a reference with directly applicable results.
It drew on the OEIS entry and simpler toothpick models without finding a solution for this variant.

The investigation proceeded by trial and error through sustained dialogue. The author directed
the search toward direct recurrences, relations for $a(2^k+i)$, coupled recurrences involving
auxiliary sequences, and comparisons with other OEIS sequences. When these attempts yielded
no general formula, he proposed working with generating functions. Several dozen
exchanges led to the representation given in this paper, followed by the asymptotic analysis.
Progress required repeated attempts and changes of direction, for which the author's mathematical
guidance was essential. With AI assistance, the author developed the proofs of the
generating-function and asymptotic results.

The AI also assisted with numerical simulations, \LaTeX{} formatting, and English expression.
The author reviewed the arguments and computations and takes full responsibility for the
mathematical content and the final text.

\end{document}